\documentclass[11pt]{article}
\usepackage{a4wide}

\newcommand{\ba}{\begin{array}}
\newcommand{\ea}{\end{array}}
\newcommand{\be}{\begin{equation}}
\newcommand{\ee}{\end{equation}}
\newcommand{\la}{\label}
\newcommand{\bea}{\begin{eqnarray}}
\newcommand{\eea}{\end{eqnarray}}

\renewcommand{\l}{\left}
\renewcommand{\r}{\right}
\newcommand{\n}{\nonumber}
\newcommand{\nn}{\nonumber \\}
\newcommand{\ds}{\displaystyle}
\newcommand{\ndots}{n=0,1,2,\ldots}
\renewcommand{\a}{\alpha}
\renewcommand{\b}{\beta}
\newcommand{\G}{\Gamma}

\newcommand{\KL}{L_n^{\a,M}(x)}
\newcommand{\GL}{L_n^{\a,M,N}(x)}
\renewcommand{\P}{P_n^{(\a,\b)}(x)}
\newcommand{\Q}{Q_n^{(\a,\b)}(x)}
\newcommand{\R}{R_n^{(\a,\b)}(x)}
\renewcommand{\S}{S_n^{(\a,\b)}(x)}
\newcommand{\GP}{P_n^{\a,\b,M,N}(x)}

\newcommand{\SGP}{P_n^{\a,\a,M,M}(x)}
\newcommand{\set}[1]{\left\{#1\right\}_{n=0}^{\infty}}
\newcommand{\hyp}[5]{\mbox{}_{#1}F_{#2}
\left(\left.{#3 \atop #4}\right|#5\right)}

\begin{document}

\title{Differential equations for generalized Jacobi polynomials}
\author{J.~Koekoek \and R.~Koekoek}
\date{}
\maketitle

\begin{abstract}
We look for differential equations of the form
\bea & &M\sum_{i=0}^{\infty}a_i(x)y^{(i)}(x)+
N\sum_{i=0}^{\infty}b_i(x)y^{(i)}(x)+
MN\sum_{i=0}^{\infty}c_i(x)y^{(i)}(x)\nn
& &{}\hspace{1cm}{}+(1-x^2)y''(x)+\l[\b-\a-(\a+\b+2)x\r]y'(x)
+n(n+\a+\b+1)y(x)=0\n\eea
satisfied by the generalized Jacobi polynomials $\set{\GP}$ which are
orthogonal on the interval $[-1,1]$ with respect to the weight function
$$\frac{\G(\a+\b+2)}{2^{\a+\b+1}\G(\a+1)\G(\b+1)}(1-x)^{\a}(1+x)^{\b}+
M\delta(x+1)+N\delta(x-1),$$
where $\a>-1$, $\b>-1$, $M\ge 0$ and $N\ge 0$. We give explicit
representations for the coefficients $\l\{a_i(x)\r\}_{i=0}^{\infty}$,
$\l\{b_i(x)\r\}_{i=0}^{\infty}$ and $\l\{c_i(x)\r\}_{i=0}^{\infty}$ and we
show that this differential equation is uniquely determined. For $M^2+N^2>0$
the order of this differential equation is infinite, except for
$\a\in\{0,1,2,\ldots\}$ or $\b\in\{0,1,2,\ldots\}$. Moreover, the order equals
$$\l\{\begin{array}{ll}
2\b+4 & \textrm{ if }\;M>0,\;N=0\;\textrm{ and }\;\b\in\{0,1,2,\ldots\}\\
2\a+4 & \textrm{ if }\;M=0,\;N>0\;\textrm{ and }\;\a\in\{0,1,2,\ldots\}\\
2\a+2\b+6 & \textrm{ if }\;M>0,\;N>0\;\textrm{ and }\;
\a,\b\in\{0,1,2,\ldots\}.\end{array}\r.$$
\end{abstract}

\vfill

\begin{tabular}{ll}
Keywords : & Differential equations, Generalized Jacobi polynomials.\\[5mm]
\multicolumn{2}{l}{1991 Mathematics Subject Classification :
Primary 33C45 ; Secondary 34A35.}
\end{tabular}

\newpage

\section{Introduction}

In \cite{Koorn} T.H.~Koornwinder introduced the polynomials $\set{\GP}$
which are orthogonal on the interval $[-1,1]$ with respect to the weight
function
$$\frac{\G(\a+\b+2)}{2^{\a+\b+1}\G(\a+1)\G(\b+1)}(1-x)^{\a}(1+x)^{\b}+
M\delta(x+1)+N\delta(x-1),$$
where $\a>-1$, $\b>-1$, $M\ge 0$ and $N\ge 0$. We call these polynomials the
generalized Jacobi polynomials, but sometimes they are also referred to as
the Jacobi-type polynomials. As a limit case he also found the generalized
Laguerre (or Laguerre-type) polynomials $\set{\KL}$ which are orthogonal on
the interval $[0,\infty)$ with respect to the weight function
$$\frac{1}{\G(\a+1)}x^{\a}e^{-x}+M\delta(x),$$
where $\a>-1$ and $M\ge 0$. These generalized Jacobi polynomials and
generalized Laguerre polynomials are related by the limit
\be\la{lim}\KL=
\lim_{\b\rightarrow\infty}P_n^{\a,\b,0,M}\l(1-\frac{2x}{\b}\r).\ee

In \cite{Dvlag} we proved that for $M>0$ the generalized Laguerre polynomials
satisfy a unique differential equation of the form
$$M\sum_{i=0}^{\infty}a_i(x)y^{(i)}(x)+xy''(x)+(\a+1-x)y'(x)+ny(x)=0,$$
where $\l\{a_i(x)\r\}_{i=0}^{\infty}$ are continuous functions on the real
line and $\l\{a_i(x)\r\}_{i=1}^{\infty}$ are independent of the degree $n$.
In \cite{Bav1} H.~Bavinck found a new method to obtain the main result of
\cite{Dvlag}. This inversion method was found in a similar way as was done in
\cite{Char} in the case of generalizations of the Charlier polynomials. See
also \cite{Invjac} for more details. In \cite{Soblag} we used this inversion
method to find all differential equations of the form
\bea & &M\sum_{i=0}^{\infty}a_i(x)y^{(i)}(x)+
N\sum_{i=0}^{\infty}b_i(x)y^{(i)}(x)\nn
& &\hspace{1cm}{}+MN\sum_{i=0}^{\infty}c_i(x)y^{(i)}(x)+
xy''(x)+(\a+1-x)y'(x)+ny(x)=0,\n\eea
where the coefficients $\l\{a_i(x)\r\}_{i=1}^{\infty}$,
$\l\{b_i(x)\r\}_{i=1}^{\infty}$ and $\l\{c_i(x)\r\}_{i=1}^{\infty}$ are
independent of $n$ and the coefficients $a_0(x)$, $b_0(x)$ and $c_0(x)$ are
independent of $x$, satisfied by the Sobolev-type Laguerre polynomials
$\set{\GL}$ which are orthogonal with respect to the inner product
$$<f,g>\;=\frac{1}{\G(\a+1)}\int_0^{\infty}x^{\a}e^{-x}f(x)g(x)dx+Mf(0)g(0)+
Nf'(0)g'(0),$$
where $\a>-1$, $M\ge 0$ and $N\ge 0$. These Sobolev-type Laguerre
polynomials $\set{\GL}$ are generalizations of the generalized Laguerre
polynomials $\set{\KL}$. In fact we have
$$L_n^{\a,M,0}(x)=\KL.$$

In this paper we will use the inversion formula found in \cite{Invjac}
to find differential equations of the form
\bea\la{DV}& &M\sum_{i=0}^{\infty}a_i(x)y^{(i)}(x)+
N\sum_{i=0}^{\infty}b_i(x)y^{(i)}(x)+
MN\sum_{i=0}^{\infty}c_i(x)y^{(i)}(x)\nn
& &{}\hspace{1cm}{}+(1-x^2)y''(x)+\l[\b-\a-(\a+\b+2)x\r]y'(x)
+n(n+\a+\b+1)y(x)=0,\eea
where the coefficients $\l\{a_i(x)\r\}_{i=1}^{\infty}$,
$\l\{b_i(x)\r\}_{i=1}^{\infty}$ and $\l\{c_i(x)\r\}_{i=1}^{\infty}$ are
independent of $n$ and the coefficients $a_0(x)$, $b_0(x)$ and $c_0(x)$
are independent of $x$, satisfied by the generalized Jacobi polynomials
$\set{\GP}$.

For $\a=\b=0$, $M>0$ and $N>0$ the generalized Jacobi polynomials
reduce to the Krall polynomials studied by L.L.~Littlejohn in
\cite{Littlejohn}. These Krall polynomials are generalizations of
the Legendre type polynomials ($\a=\b=0$ and $N=M>0$) found by
H.L.~Krall in \cite{H.L.Krall1} and \cite{H.L.Krall2}. See also
\cite{A.M.Krall}. In \cite{Littlejohn} it is shown that the Krall
polynomials satisfy a sixth order differential equation of the
form (\ref{DV}). For $\a>-1$, $\b=0$, $M>0$ and $N=0$ or for
$\a=0$, $\b>-1$, $M=0$ and $N>0$ the generalized Jacobi polynomials
reduce to the Jacobi type polynomials which satisfy a fourth order
differential equation of the form (\ref{DV})~; see also
\cite{A.M.Krall}, \cite{H.L.Krall1} and \cite{H.L.Krall2}.

We emphasize that the case $\b=\a$ and $N=M$ is special in the
sense that we can also find differential equations of the form
\be\la{DVSym}M\sum_{i=0}^{\infty}d_i(x)y^{(i)}(x)+(1-x^2)y''(x)
-2(\a+1)xy'(x)+n(n+2\a+1)y(x)=0,\ee
where the coefficients $\l\{d_i(x)\r\}_{i=1}^{\infty}$ are
independent of $n$ and $d_0(x)$ is independent of $x$,
satisfied by the symmetric generalized ultraspherical polynomials
$\set{\SGP}$. The Legendre type polynomials for instance satisfy
a fourth order differential equation of the form (\ref{DVSym}).
See \cite{A.M.Krall}, \cite{H.L.Krall1} and \cite{H.L.Krall2}. In
\cite{Symjac} we found all differential equations of the form
(\ref{DVSym}) satisfied by the polynomials $\set{\SGP}$ for $\a>-1$
and $M\ge 0$. In \cite{Madrid} we applied the special case $\b=\a$
of the Jacobi inversion formula to solve the systems of equations
obtained in \cite{Symjac}.

\section{The main results}

We look for all differential equations of the form (\ref{DV}) satisfied by
the generalized Jacobi polynomials $\set{\GP}$. A representation of these
orthogonal polynomials will be given in section 5. We emphasize that we
demand that the coefficients $\l\{a_i(x)\r\}_{i=1}^{\infty}$,
$\l\{b_i(x)\r\}_{i=1}^{\infty}$ and $\l\{c_i(x)\r\}_{i=1}^{\infty}$ are
independent of the degree $n$ and that $a_0(x)$, $b_0(x)$ and $c_0(x)$ do
not depend on $x$. Therefore we will use the following notations~:
$$a_0(x)=a_0(n,\a,\b),\quad b_0(x)=b_0(n,\a,\b),\quad c_0(x)=c_0(n,\a,\b),
\;\ndots$$
and
$$a_i(x)=a_i(\a,\b,x),\quad b_i(x)=b_i(\a,\b,x),\quad c_i(x)=c_i(\a,\b,x),
\;i=1,2,3,\ldots.$$

We will apply a general theorem by H.~Bavinck to prove that for $\a>-1$,
$\b>-1$, $M\ge 0$ and $N\ge 0$ the polynomials $\set{\GP}$ satisfy a
unique differential equation of the form (\ref{DV}), where
\be\la{anul}a_0(0,\a,\b)=0,\quad a_0(n,\a,\b)=(\a+\b+2)
\frac{(\b+3)_{n-1}(\a+\b+3)_{n-1}}{(\a+1)_{n-1}\,(n-1)!},\;n=1,2,3,\ldots,\ee
\be\la{bnul}b_0(0,\a,\b)=0,\quad b_0(n,\a,\b)=(\a+\b+2)
\frac{(\a+3)_{n-1}(\a+\b+3)_{n-1}}{(\b+1)_{n-1}\,(n-1)!},\;n=1,2,3,\ldots,\ee
\bea\la{cnul}& &c_0(0,\a,\b)=c_0(1,\a,\b)=0\;\textrm{ and}\nn
& &c_0(n,\a,\b)=\frac{(\a+\b+2)^2(\a+\b+3)}{(\a+1)(\b+1)}\nn
& &{}\hspace{3cm}\times\frac{(\a+\b+4)_{n-1}}{(n-1)!}
\frac{(\a+\b+4)_{n-2}}{(n-2)!},\;n=2,3,4,\ldots.\eea

Further we will show that
\bea\la{a} & &a_i(\a,\b,x)=-(\a+\b+2)2^i
\sum_{\ell=0}^{i-1}(-1)^{\ell}\frac{(\b+3)_{i-\ell-1}(-\b-2)_{i-\ell-1}}
{(\a+1)_{i-\ell-1}\,(i-\ell)!\,(i-\ell-1)!\,\ell!}\nn
& &{}\hspace{3cm}\times
\hyp{3}{2}{-\ell,\a+\b+3,\b+i-\ell+2}{\a+i-\ell,i-\ell+1}{1}
\l(\frac{x+1}{2}\r)^{\ell+1}\eea
and
\bea\la{b} & &b_i(\a,\b,x)=(\a+\b+2)(-2)^i
\sum_{\ell=0}^{i-1}\frac{(\a+3)_{i-\ell-1}(-\a-2)_{i-\ell-1}}
{(\b+1)_{i-\ell-1}\,(i-\ell)!\,(i-\ell-1)!\,\ell!}\nn
& &{}\hspace{3cm}\times
\hyp{3}{2}{-\ell,\a+\b+3,\a+i-\ell+2}{\b+i-\ell,i-\ell+1}{1}
\l(\frac{x-1}{2}\r)^{\ell+1}\eea
for $i=1,2,3,\ldots$ and that
\be\la{c}c_1(\a,\b,x)=0\;\textrm{ and }\;c_i(\a,\b,x)=
c_i^{(1)}(\a,\b,x)+c_i^{(2)}(\a,\b,x),\;i=2,3,4,\ldots,\ee
where for $i=2,3,4,\ldots$
\bea\la{c1} & &c_i^{(1)}(\a,\b,x)=
-\frac{(\a+\b+2)^2(\a+\b+3)(\a+\b+4)}{(\a+1)(\b+1)i}(x^2-1)2^{i-2}\nn
& &{}\hspace{2cm}\times\sum_{\ell=0}^{i-2}(-1)^{\ell}
\frac{(\b+3)_{i-\ell-2}(-\a-\b-3)_{i-\ell-2}}
{(i-\ell-1)!\,(i-\ell-2)!\,\ell!\,(i-\ell-1)!}\nn
& &{}\hspace{3cm}\times\hyp{4}{3}{-\ell,\a+\b+5,\a+\b+4,\b+i-\ell+1}
{\b+3,i-\ell,i-\ell}{1}\l(\frac{x+1}{2}\r)^{\ell+1}\eea
and
\bea\la{c2} & &c_i^{(2)}(\a,\b,x)=
\frac{(\a+\b+2)^2(\a+\b+3)(\a+\b+4)}{(\a+1)(\b+1)i}(x^2-1)(-2)^{i-2}\nn
& &{}\hspace{2cm}\times\sum_{\ell=0}^{i-2}
\frac{(\a+3)_{i-\ell-2}(-\a-\b-3)_{i-\ell-2}}
{(i-\ell-1)!\,(i-\ell-2)!\,\ell!\,(i-\ell-1)!}\nn
& &{}\hspace{3cm}\times\hyp{4}{3}{-\ell,\a+\b+5,\a+\b+4,\a+i-\ell+1}
{\a+3,i-\ell,i-\ell}{1}\l(\frac{x-1}{2}\r)^{\ell+1}.\eea

Note that we have
\be\la{symab}a_i(\a,\b,x)=(-1)^ib_i(\b,\a,-x),\;i=1,2,3,\ldots\ee
and
\be\la{symc}c_i^{(1)}(\a,\b,x)=(-1)^ic_i^{(2)}(\b,\a,-x),\;i=2,3,4,\ldots.\ee

Finally we will show that for $\a>-1$, $\b>-1$ and $M^2+N^2>0$ the order of
the differential equation (\ref{DV}) will be infinite in general. Only for
nonnegative integer values of $\a$ or $\b$ finite order can occur. Moreover,
the order of the differential equation equals
$$\l\{\begin{array}{ll}
2\b+4 & \textrm{ if }\;M>0,\;N=0\;\textrm{ and }\;\b\in\{0,1,2,\ldots\}\\
2\a+4 & \textrm{ if }\;M=0,\;N>0\;\textrm{ and }\;\a\in\{0,1,2,\ldots\}\\
2\a+2\b+6 & \textrm{ if }\;M>0,\;N>0\;\textrm { and }\;
\a,\b\in\{0,1,2,\ldots\}.\end{array}\r.$$

In fact, we will show that
\be\la{afbra}a_i(\a,\b,x)=0,\;i>2\b+4\;\textrm{ if }\;\b\in\{0,1,2,\ldots\},\ee
\be\la{afbrb}b_i(\a,\b,x)=0,\;i>2\a+4\;\textrm{ if }\;\a\in\{0,1,2,\ldots\}\ee
and
\be\la{afbrc}c_i(\a,\b,x)=0,\;i>2\a+2\b+6\;\textrm{ if }\;\a,\b\in\{0,1,2,\ldots\}.\ee
Further we have
\be\la{akop}a_{2\b+4}(\a,\b,x)=-\frac{1}{(\a+1)_{\b+1}}
\frac{(x^2-1)^{\b+2}}{(\b+2)!},\;\b\in\{0,1,2,\ldots\},\ee
\be\la{bkop}b_{2\a+4}(\a,\b,x)=-\frac{1}{(\b+1)_{\a+1}}
\frac{(x^2-1)^{\a+2}}{(\a+2)!},\;\a\in\{0,1,2,\ldots\}\ee
and
\be\la{ckop}c_{2\a+2\b+6}(\a,\b,x)=-\frac{\a+\b+2}{(\a+1)(\b+1)}
\frac{(x^2-1)^{\a+\b+3}}{(\a+\b+1)!\,(\a+\b+3)!},
\;\a,\b\in\{0,1,2,\ldots\}.\ee

\section{The classical Jacobi polynomials}

In this section we list the definitions and some properties of the classical
Jacobi polynomials which we will use in this paper. For details the reader
is referred to \cite{Chihara}, \cite{Askey} and \cite{Szego}.

The classical Jacobi polynomials $\set{\P}$ can be defined by
\bea\la{defJac1}\P&=&\sum_{k=0}^n\frac{(n+\a+\b+1)_k}{k!}
\frac{(\a+k+1)_{n-k}}{(n-k)!}\l(\frac{x-1}{2}\r)^k,\;\ndots\\
&=&\la{defJac2}(-1)^n\sum_{k=0}^n\frac{(-n-k-\a-\b)_k}{k!}
\frac{(-n-\a)_{n-k}}{(n-k)!}\l(\frac{x-1}{2}\r)^k,\;\ndots\eea
for all $\a$ and $\b$. The Jacobi polynomials satisfy the symmetry relation
\be\la{symJac}P_n^{(\a,\b)}(x)=(-1)^nP_n^{(\b,\a)}(-x),\;\ndots.\ee
From (\ref{defJac1}) and (\ref{symJac}) we easily find for $\ndots$
\be\la{randen}P_n^{(\a,\b)}(1)=\frac{(\a+1)_n}{n!}\;\textrm{ and }\;
P_n^{(\a,\b)}(-1)=(-1)^n\frac{(\b+1)_n}{n!}\ee
and
\be\la{diff}D^i\P=\frac{(n+\a+\b+1)_i}{2^i}P_{n-i}^{(\a+i,\b+i)}(x),
\;i=0,1,2,\ldots,n,\ee
where $D=\ds\frac{d}{dx}$ denotes the differentiation operator. These Jacobi
polynomials satisfy the linear second order differential equation
\be\la{dvJac}(1-x^2)y''(x)+\l[\b-\a-(\a+\b+2)x\r]y'(x)
+n(n+\a+\b+1)y(x)=0.\ee

By using the definition (\ref{defJac1}) and the symmetry relation (\ref{symJac})
it is not very difficult to derive the following relations
\be\la{rel1}P_n^{(\a+1,\b)}(x)-P_n^{(\a,\b+1)}(x)=P_{n-1}^{(\a+1,\b+1)}(x),
\;n=1,2,3,\ldots,\ee
\be\la{rel2a}n\P-(n+\a)P_{n-1}^{(\a,\b+1)}(x)=(x-1)D\P,\;n=1,2,3,\ldots,\ee
\be\la{rel2b}n\P+(n+\b)P_{n-1}^{(\a+1,\b)}(x)=(x+1)D\P,\;n=1,2,3,\ldots,\ee
\be\la{rel3a}(n+\a+1)\P-(\a+1)P_n^{(\a+1,\b)}(x)=
(n+\b)\l(\frac{x-1}{2}\r)P_{n-1}^{(\a+2,\b)}(x),\;n=1,2,3,\ldots\ee
and
\be\la{rel3b}(n+\b+1)\P-(\b+1)P_n^{(\a,\b+1)}(x)=
(n+\a)\l(\frac{x+1}{2}\r)P_{n-1}^{(\a,\b+2)}(x),\;n=1,2,3,\ldots.\ee
Note that the differential equation (\ref{dvJac}) implies that
\bea\la{rel4}& &n(n+\a+\b+1)\P-\l[(\b+1)(x-1)+(\a+1)(x+1)\r]D\P\nn
& &{}\hspace{1cm}=(x^2-1)D^2\P,\;\ndots.\eea
By using Leibniz' rule we also have for $\ndots$ and $i=0,1,2,\ldots$
\bea\la{rel5}& &(1-x^2)D^{i+2}\P+\l[\b-\a-(\a+\b+2i+2)x\r]D^{i+1}\P\nn
& &{}\hspace{1cm}{}+(n-i)(n+\a+\b+i+1)D^i\P=0.\eea

\section{Some inversion, summation and transformation formulas}

In this section we will give some inversion formulas which we will need in
this paper. Further we derive some summation formulas which we will use.
Finally we give two transformation formulas which will be used in section 8
of this paper.

Let $\a>-1$ and $\b>-1$.

In this paper we have to deal with systems of equations of the form
\be\la{inv1}\sum_{i=1}^{\infty}A_i(x)D^i\P=F_n(x),\;n=1,2,3,\ldots,\ee
where the coefficients $\l\{A_i(x)\r\}_{i=1}^{\infty}$ are independent of
$n$. In \cite{Invjac} we have shown that this system of equations has a
unique solution given by
\be\la{opl1}A_i(x)=2^i\sum_{j=1}^i\frac{\a+\b+2j+1}{(\a+\b+j+1)_{i+1}}
P_{i-j}^{(-\a-i-1,-\b-i-1)}(x)F_j(x),\;i=1,2,3,\ldots.\ee

We will also need a variant of this inversion formula. In a similar way we
may also conclude that a system of equations of the form
\be\la{inv2}\sum_{i=0}^{\infty}B_i(x)D^i\P=G_n(x),\;\ndots,\ee
where the coefficients $\l\{B_i(x)\r\}_{i=0}^{\infty}$ are independent of
$n$ has a unique solution given by
\be\la{opl2}B_i(x)=2^i\sum_{j=0}^i\frac{\a+\b+2j+1}{(\a+\b+j+1)_{i+1}}
P_{i-j}^{(-\a-i-1,-\b-i-1)}(x)G_j(x),\;i=0,1,2,\ldots.\ee
The case $\a+\b+1=0$ must be understood by continuity.

Let $N$ denote a positive integer. Now we consider the $(N\times N)$-matrix
$A$ defined by
\be\la{A1}A=(a_{ij})_{i,j=1}^N\;\textrm{ with }\;
a_{ij}=\l\{\begin{array}{ll}i, & i=j\\z, & i=j+1\\0, & \textrm{otherwise.}
\end{array}\r.\ee
Since $\textrm{det}(A)=N!\ne 0$ this matrix is invertible for every $z$. We
will show that its inverse is given by
\be\la{B1}A^{-1}=B=(b_{ij})_{i,j=1}^N\;\textrm{ with }\;
b_{ij}=\l\{\begin{array}{ll}
\ds\frac{(-1)^{i-j}(j-1)!\,z^{i-j}}{i!}, & i\ge j\\
0, & i<j.\end{array}\r.\ee
To prove this we write
$$AB=C=(c_{ij})_{i,j=1}^N\;\textrm{ with }\;
c_{ij}=\sum_{k=1}^Na_{ik}b_{kj}$$
and we will show that $C=I$, the identity matrix. For $N=1$ this is
trivially true. For $N\in\{2,3,4,\ldots\}$ we have
$$c_{1j}=a_{11}b_{1j}=b_{1j}\;\textrm{ and }\;
c_{ij}=a_{i,i-1}b_{i-1,j}+a_{ii}b_{ij}=zb_{i-1,j}+ib_{ij}\;\textrm{ for }\;i\ge 2.$$
Hence $c_{ij}=0$ if $i<j$,
$$c_{11}=b_{11}=1,\;c_{ii}=0+i\frac{(i-1)!}{i!}=1\;\textrm{ for }i\ge 2$$
and
$$c_{ij}=z\frac{(-1)^{i-j-1}(j-1)!\,z^{i-j-1}}{(i-1)!}
+i\frac{(-1)^{i-j}(j-1)!\,z^{i-j}}{i!}=0\;\textrm{ for }\;
i\ge j+1.$$
This proves (\ref{B1}).

We also need the following matrix inverse. Let $N$ denote a positive integer
again and consider the $(N\times N)$-matrix $A$ defined by
\be\la{A2}A=(a_{ij})_{i,j=1}^N\;\textrm{ with }\;
a_{ij}=\l\{\begin{array}{ll}i(i+1), & i=j\\2ix, & i=j+1\\x^2-1, & i=j+2\\
0, & \textrm{otherwise.}\end{array}\r.\ee
Since $\textrm{det}(A)=N!\,(N+1)!\ne 0$ this matrix is invertible for every
$x$. We will show that its inverse is given by
\bea\la{B2}& &A^{-1}=B=(b_{ij})_{i,j=1}^N\nn
& &{}\hspace{1cm}\textrm{with }\;b_{ij}=\l\{\begin{array}{ll}
\ds\frac{(-1)^{i-j}(j-1)!\,\l[(x+1)^{i-j+1}-(x-1)^{i-j+1}\r]}{2(i+1)!}, &
i\ge j\\0, & i<j.\end{array}\r.\eea
To prove this we write again
$$AB=C=(c_{ij})_{i,j=1}^N\;\textrm{ with }\;
c_{ij}=\sum_{k=1}^Na_{ik}b_{kj}$$
and again we will show that $C=I$, the identity matrix. For $N=1$ this
is trivially true and for $N=2$ we find that
$$C=AB=\l(\begin{array}{cc}2 & 0\\4x & 6\end{array}\r)
\l(\begin{array}{cc}\frac{1}{2} & 0\\-\frac{1}{3}x & \frac{1}{6}\end{array}\r)
=\l(\begin{array}{cc}1 & 0\\0 & 1\end{array}\r)=I.$$
For $N\in\{3,4,5,\ldots\}$ we have
$$c_{1j}=a_{11}b_{1j}=2b_{1j},\;
c_{2j}=a_{21}b_{1j}+a_{22}b_{2j}=4xb_{1j}+6b_{2j}$$
and
$$c_{ij}=a_{i,i-2}b_{i-2,j}+a_{i,i-1}b_{i-1,j}+a_{ii}b_{ij}
=(x^2-1)b_{i-2,j}+2ixb_{i-1,j}+i(i+1)b_{ij}\;\textrm{ for }\;i\ge 3.$$
Hence $c_{ij}=0$ if $i<j$,
$$c_{11}=2b_{11}=1,\;c_{22}=0+6b_{22}=1,$$
$$c_{ii}=0+0+i(i+1)\frac{(i-1)!\,2}{2(i+1)!}=1\;\textrm{ for }\;i\ge 3,$$
$$c_{21}=4xb_{11}+6b_{21}=0,\;
c_{i,i-1}=0+2ix\frac{(i-2)!\,2}{2\,i!}-i(i+1)\frac{(i-2)!\,4x}{2(i+1)!}=0
\;\textrm{ for }\;i\ge 3$$
and by using $2x=(x+1)+(x-1)$
\bea c_{ij}&=&(x^2-1)
\frac{(-1)^{i-j}(j-1)!\,\l[(x+1)^{i-j-1}-(x-1)^{i-j-1}\r]}{2(i-1)!}\nn
& &{}\hspace{1cm}+2ix
\frac{(-1)^{i-j-1}(j-1)!\,\l[(x+1)^{i-j}-(x-1)^{i-j}\r]}{2\,i!}\nn
& &{}\hspace{2cm}+i(i+1)
\frac{(-1)^{i-j}(j-1)!\,\l[(x+1)^{i-j+1}-(x-1)^{i-j+1}\r]}{2\,(i+1)!}\nn
&=&\frac{(-1)^{i-j}(j-1)!}{2(i-1)!}
\l[(x-1)(x+1)^{i-j}-(x+1)(x-1)^{i-j}-(x+1)^{i-j+1}\r.\nn
& &{}\hspace{4cm}+(x+1)(x-1)^{i-j}-(x-1)(x+1)^{i-j}+(x-1)^{i-j+1}\nn
& &{}\hspace{5cm}\l.{}+(x+1)^{i-j+1}-(x-1)^{i-j+1}\r]\nn
&=&0\;\textrm{ for }\;i\ge j+2.\n\eea
This proves (\ref{B2}).

We will also need the well-known Vandermonde summation formula
\be\la{Van}\hyp{2}{1}{-n,a}{b}{1}=\frac{(b-a)_n}{(b)_n},\;(b)_n\ne 0,
\;\ndots,\ee
which can be found in \cite{Bailey} and \cite{Slater} for instance. We also
need the following summation formulas~:
\bea\la{F}F_n(a,b)&=&\sum_{k=0}^n\frac{(a)_k(b)_k}{(b-a+1)_k\,k!}(b+2k)\\
&=&\la{sum1}\frac{(a+1)_n(b)_{n+1}}{(b-a+1)_n\,n!},\;\ndots\eea
and
\bea\la{sum2}& &\sum_{k=0}^n\frac{(-n)_k(a)_k(b)_k(c)_k}
{(b+n+1)_k(b-a+1)_k(b-c+1)_k\,k!}(b+2k)\nn
& &{}\hspace{1cm}=\frac{(b)_{n+1}(b-a-c+1)_n}{(b-a+1)_n(b-c+1)_n},
\;\ndots.\eea
Formula (\ref{sum1}) can easily be proved by using mathematical
induction. Formula (\ref{sum2}) can be proved by using the well-known
summation formula for a terminating well-poised ${}_5F_4$~:
$$\hyp{5}{4}{-n,a,b,c,\frac{1}{2}b+1}{b+n+1,b-a+1,b-c+1,\frac{1}{2}b}{1}=
\frac{(b+1)_n(b-a-c+1)_n}{(b-a+1)_n(b-c+1)_n},\;\ndots.$$
This formula can be found in \cite{Bailey} and \cite{Slater} for instance.
Note that (\ref{sum1}) follows from (\ref{sum2}) by setting
$c=b+n+1$.

Finally we will need the following transformation formula (see for instance
\cite{Luke}, section 9.1, formula (34))
\bea\la{trans1}& &\hyp{4}{3}{a,b,c,p}{d,e,q}{z}=\sum_{n=0}^{\infty}
(-1)^n\frac{(a)_n(b)_n(c)_n(q-p)_n}{(d)_n(e)_n(q)_n\,n!}z^n\nn
& &{}\hspace{5cm}\times\hyp{3}{2}{n+a,n+b,n+c}{n+d,n+e}{z},
\;\textrm{Re}(z)<\frac{1}{2}.\eea
As a special case we also have
\be\la{trans2}\hyp{3}{2}{a,b,p}{c,q}{z}=\sum_{n=0}^{\infty}
(-1)^n\frac{(a)_n(b)_n(q-p)_n}{(c)_n(q)_n\,n!}z^n\,\hyp{2}{1}{n+a,n+b}{n+c}{z},
\;\textrm{Re}(z)<\frac{1}{2}.\ee

\section{The generalized Jacobi polynomials}

Let $\a>-1$, $\b>-1$, $M\ge 0$ and $N\ge 0$. In \cite{Koorn} it is shown
that the generalized Jacobi polynomials $\set{\GP}$ can be written as
\be\la{def}\GP=\P+M\Q+N\R+MN\S,\;\ndots,\ee
where
$$Q_0^{(\a,\b)}(x)=R_0^{(\a,\b)}(x)=S_0^{(\a,\b)}(x)=0$$
and for $n=1,2,3,\ldots$
\bea\la{defQ}& &\Q=\frac{(\b+2)_{n-1}(\a+\b+2)_{n-1}}{(\a+1)_n\,n!}\nn
& &{}\hspace{3cm}\times\l[n(n+\a+\b+1)\P-(\b+1)(x-1)D\P\r],\eea
\bea\la{defR}& &\R=\frac{(\a+2)_{n-1}(\a+\b+2)_{n-1}}{(\b+1)_n\,n!}\nn
& &{}\hspace{3cm}\times\l[n(n+\a+\b+1)\P-(\a+1)(x+1)D\P\r]\eea
and
\bea\la{defS}& &\S=\frac{1}{(\a+1)(\b+1)}
\frac{(\a+\b+2)_n(\a+\b+2)_{n-1}}{n!\,(n-1)!}\nn
& &{}\hspace{3cm}\times
\l[n(n+\a+\b+1)\P\r.\nn
& &{}\hspace{5cm}\l.{}-\l\{(\b+1)(x-1)+(\a+1)(x+1)\r\}D\P\r].\eea

First of all we remark that the generalized Jacobi polynomials satisfy
the symmetry relation (see \cite{Koorn})
\be\la{sym}P_n^{\a,\b,M,N}(x)=(-1)^nP_n^{\b,\a,N,M}(-x),\;\ndots,\ee
which implies that
\be\la{symQR}\Q=(-1)^nR_n^{(\b,\a)}(-x),\;\ndots\ee
and
$$\S=(-1)^nS_n^{(\b,\a)}(-x),\;\ndots.$$

From (\ref{defQ}) and (\ref{defR}) it follows that
\be\la{Qplus}Q_n^{(\a,\b)}(1)=
\frac{(\b+2)_{n-1}(\a+\b+2)_n}{(\a+1)_n\,(n-1)!}
P_n^{(\a,\b)}(1),\;n=1,2,3,\ldots\ee
and
\be\la{Rmin}R_n^{(\a,\b)}(-1)=
\frac{(\a+2)_{n-1}(\a+\b+2)_n}{(\b+1)_n\,(n-1)!}
P_n^{(\a,\b)}(-1),\;n=1,2,3,\ldots.\ee
These two formulas will be used in the next section.

Now we use (\ref{rel2a}), (\ref{rel1}) and (\ref{rel2b}) to obtain
for $n=1,2,3,\ldots$
\bea & &n(n+\a+\b+1)\P-(\b+1)(x-1)D\P\nn
&=&(n+\a)\l[n\P+(\b+1)P_{n-1}^{(\a,\b+1)}(x)\r]\nn
&=&(n+\a)\l[nP_n^{(\a-1,\b+1)}(x)+(n+\b+1)P_{n-1}^{(\a,\b+1)}(x)\r]\nn
&=&(n+\a)(x+1)DP_n^{(\a-1,\b+1)}(x).\n\eea
Hence from (\ref{defQ}) we obtain the following representations
for $\Q$~:
\bea\Q&=&\la{defQ2}\frac{(\b+2)_{n-1}(\a+\b+2)_{n-1}}{(\a+1)_{n-1}\,n!}
\l[n\P+(\b+1)P_{n-1}^{(\a,\b+1)}(x)\r]\\
&=&\la{defQ1}\frac{(\b+2)_{n-1}(\a+\b+2)_{n-1}}{(\a+1)_{n-1}\,n!}
(x+1)DP_n^{(\a-1,\b+1)}(x)\eea
for $n=1,2,3,\ldots$. In a similar way from (\ref{defR}) or by
using the symmetry relation (\ref{symQR}) we find the following
representations for $\R$~:
\bea\R&=&\la{defR2}\frac{(\a+2)_{n-1}(\a+\b+2)_{n-1}}{(\b+1)_{n-1}\,n!}
\l[n\P-(\a+1)P_{n-1}^{(\a+1,\b)}(x)\r]\\
&=&\la{defR1}\frac{(\a+2)_{n-1}(\a+\b+2)_{n-1}}{(\b+1)_{n-1}\,n!}
(x-1)DP_n^{(\a+1,\b-1)}(x)\eea
for $n=1,2,3,\ldots$.

And if we use (\ref{rel4}) we easily find from (\ref{defS}) that
\bea\la{defS1}& &\S=\frac{1}{(\a+1)(\b+1)}\nn
& &{}\hspace{2cm}\times\frac{(\a+\b+2)_n(\a+\b+2)_{n-1}}{n!\,(n-1)!}
(x^2-1)D^2\P,\;n=1,2,3,\ldots.\eea

Note that the representations (\ref{defQ2}) and (\ref{defR2}) imply that
for $n=1,2,3,\ldots$ we have
\be\la{q}\Q=\sum_{k=0}^nq_{n,k}^{(\a,\b)}P_k^{(\a,\b)}(x)\;\textrm{ with }\;
q_{n,n}^{(\a,\b)}=\frac{(\b+2)_{n-1}(\a+\b+2)_{n-1}}{(\a+1)_{n-1}\,(n-1)!}\ee
and
\be\la{r}\R=\sum_{k=0}^nr_{n,k}^{(\a,\b)}P_k^{(\a,\b)}(x)\;\textrm{ with }\;
r_{n,n}^{(\a,\b)}=\frac{(\a+2)_{n-1}(\a+\b+2)_{n-1}}{(\b+1)_{n-1}\,(n-1)!}.\ee
By using (\ref{rel2a}) and (\ref{rel2b}) we also find from (\ref{defS})
that for $n=1,2,3,\ldots$ we have
$$\S=\sum_{k=0}^ns_{n,k}^{(\a,\b)}P_k^{(\a,\b)}(x)$$
with
\be\la{s}s_{n,n}^{(\a,\b)}=\frac{n(n-1)}{(\a+1)(\b+1)}
\frac{(\a+\b+2)_n(\a+\b+2)_{n-1}}{n!\,(n-1)!}.\ee

\section{The existence and uniqueness of the differential equation
and the 'eigenvalue' coefficients}

First of all we set
$$\lambda_n=n(n+\a+\b+1),\;\ndots,$$
which implies that $\lambda_0=0$ and
\be\la{verschil}\lambda_n-\lambda_{n-1}=2n+\a+\b,\;n=1,2,3,\ldots.\ee

In \cite{Bav3} H.~Bavinck proved a theorem concerning differential or
difference equations satisfied by certain orthogonal polynomials. This
result can be applied to the generalized Jacobi polynomials $\set{\GP}$.
In that case for $\a>-1$, $\b>-1$, $M\ge 0$ and $N\ge 0$ his result reads
as follows~:

\vspace{5mm}

{\bf Theorem} (H.~Bavinck). {\em If
\be\la{cond1}P_n^{(\a,\b)}(1)\ne 0\;\textrm{ and }\;
P_n^{(\a,\b)}(-1)\ne 0,\;\ndots\ee
and
\be\la{cond2}P_n^{(\a,\b)}(1)+MQ_n^{(\a,\b)}(1)\ne 0\;\textrm{ and }\;
P_n^{(\a,\b)}(-1)+NR_n^{(\a,\b)}(-1)\ne 0,\;\ndots\ee
then the generalized Jacobi polynomials given by (\ref{def}) satisfy a
unique differential equation of the form (\ref{DV}), where
$$a_0(0,\a,\b)=b_0(0,\a,\b)=c_0(0,\a,\b)=0,$$
$$a_0(n,\a,\b)=\sum_{j=1}^n(\lambda_j-\lambda_{j-1})
q_{j,j}^{(\a,\b)},\;n=1,2,3,\ldots,$$
$$b_0(n,\a,\b)=\sum_{j=1}^n(\lambda_j-\lambda_{j-1})
r_{j,j}^{(\a,\b)},\;n=1,2,3,\ldots$$
and
$$c_0(n,\a,\b)=\sum_{j=1}^n(\lambda_j-\lambda_{j-1})
s_{j,j}^{(\a,\b)},\;n=1,2,3,\ldots,$$
where $q_{j,j}^{(\a,\b)}$, $r_{j,j}^{(\a,\b)}$ and $s_{j,j}^{(\a,\b)}$
are given by (\ref{q}), (\ref{r}) and (\ref{s}).}

\vspace{5mm}

Since $\a>-1$, $\b>-1$, $P_0^{(\a,\b)}(x)=1$ and
$Q_0^{(\a,\b)}(x)=R_0^{(\a,\b)}(x)=0$ it easily follows from
(\ref{randen}) that condition (\ref{cond1}) is satisfied.
Since also $M\ge 0$ and $N\ge 0$ we conclude, by using (\ref{Qplus})
and (\ref{Rmin}), that condition (\ref{cond2}) is satisfied too.

By using (\ref{q}), (\ref{r}), (\ref{s}) and (\ref{verschil}) we
find that
$$a_0(n,\a,\b)=\sum_{j=1}^n\frac{(\b+2)_{j-1}(\a+\b+2)_{j-1}}
{(\a+1)_{j-1}\,(j-1)!}(2j+\a+\b),\;n=1,2,3,\ldots,$$
$$b_0(n,\a,\b)=\sum_{j=1}^n\frac{(\a+2)_{j-1}(\a+\b+2)_{j-1}}
{(\b+1)_{j-1}\,(j-1)!}(2j+\a+\b),\;n=1,2,3,\ldots,$$
$c_0(1,\a,\b)=0$ and
\bea & &c_0(n,\a,\b)=\frac{(\a+\b+2)^2}{(\a+1)(\b+1)}\nn
& &{}\hspace{3cm}\times\sum_{j=2}^n\frac{(\a+\b+3)_{j-1}}{(j-1)!}
\frac{(\a+\b+3)_{j-2}}{(j-2)!}(2j+\a+\b),\;n=2,3,4,\ldots.\n\eea

Note that
$$a_0(n,\a,\b)=\sum_{k=0}^{n-1}\frac{(\b+2)_k(\a+\b+2)_k}{(\a+1)_k\,k!}
(2k+\a+\b+2)=F_{n-1}(\b+2,\a+\b+2),\;n=1,2,3,\ldots,$$
$$b_0(n,\a,\b)=\sum_{k=0}^{n-1}\frac{(\a+2)_k(\a+\b+2)_k}{(\b+1)_k\,k!}
(2k+\a+\b+2)=F_{n-1}(\a+2,\a+\b+2),\;n=1,2,3,\ldots$$
and since $(k+1)!=(2)_k,\;k=0,1,2,\ldots$
\bea c_0(n,\a,\b)&=&\frac{(\a+\b+2)^2(\a+\b+3)}{(\a+1)(\b+1)}
\sum_{k=0}^{n-2}\frac{(\a+\b+3)_k(\a+\b+4)_k}{(2)_k\,k!}(2k+\a+\b+4)\nn
&=&\frac{(\a+\b+2)^2(\a+\b+3)}{(\a+1)(\b+1)}F_{n-2}(\a+\b+3,\a+\b+4),
\;n=2,3,4,\ldots,\n\eea
where $F_n(a,b)$ is given by (\ref{F}). Now we use the summation formula
(\ref{sum1}) to obtain (\ref{anul}), (\ref{bnul}) and (\ref{cnul}).

\section{The computation of the other coefficients}

First of all we remark that the symmetry relation (\ref{sym}) implies that
\be\la{symnul}a_0(n,\a,\b)=b_0(n,\b,\a)\;\textrm{ and }
c_0(n,\a,\b)=c_0(n,\b,\a),\;\ndots\ee
and
$$a_i(\a,\b,x)=(-1)^ib_i(\b,\a,-x)\;\textrm{ and }\;
c_i(\a,\b,x)=(-1)^ic_i(\b,\a,-x),\;i=1,2,3,\ldots.$$
Hence we have (\ref{symab}). Note that in the preceding section we
have already determined the 'eigenvalue' coefficients $a_0(n,\a,\b)$,
$b_0(n,\a,\b)$ and $c_0(n,\a,\b)$. From (\ref{anul}), (\ref{bnul})
and (\ref{cnul}) it is clear that (\ref{symnul}) is satisfied.

In order to compute the other coefficients $\l\{a_i(x)\r\}_{i=1}^{\infty}$,
$\l\{b_i(x)\r\}_{i=1}^{\infty}$ and $\l\{c_i(x)\r\}_{i=1}^{\infty}$ we set
$y(x)=\GP$ in the differential equation (\ref{DV}) and use (\ref{def}) and
the fact that the classical Jacobi polynomials satisfy the
differential equation (\ref{dvJac}) to obtain for $\ndots$
\bea\la{invul} & &M\sum_{i=0}^{\infty}a_i(x)D^i\l[\P+M\Q+N\R+MN\S\r]\nn
& &{}+N\sum_{i=0}^{\infty}b_i(x)D^i\l[\P+M\Q+N\R+MN\S\r]\nn
& &{}+MN\sum_{i=0}^{\infty}c_i(x)D^i\l[\P+M\Q+N\R+MN\S\r]\nn
& &{}+M\l[(1-x^2)D^2\Q+\l[\b-\a-(\a+\b+2)x\r]D\Q\r.\nn
& &{}\hspace{7cm}\l.{}+n(n+\a+\b+1)\Q\r]\nn
& &{}+N\l[(1-x^2)D^2\R+\l[\b-\a-(\a+\b+2)x\r]D\R\r.\nn
& &{}\hspace{7cm}\l.{}+n(n+\a+\b+1)\R\r]\nn
& &{}+MN\l[(1-x^2)D^2\S+\l[\b-\a-(\a+\b+2)x\r]D\S\r.\nn
& &{}\hspace{7cm}\l.{}+n(n+\a+\b+1)\S\r]=0.\eea
Since $P_0^{(\a,\b)}(x)=1$,
$Q_0^{(\a,\b)}(x)=R_0^{(\a,\b)}(x)=S_0^{(\a,\b)}(x)=0$ and
$a_0(0,\a,\b)=b_0(0,\a,\b)=c_0(0,\a,\b)=0$ this is trivially true for $n=0$.

Now we use (\ref{defQ2}), (\ref{defR2}) and (\ref{rel1}) to find for
$n=1,2,3,\ldots$
$$\Q=\frac{(\b+2)_{n-1}(\a+\b+2)_{n-1}}{(\a+1)_{n-1}\,n!}
\l[(n+\b+1)\P-(\b+1)P_n^{(\a-1,\b+1)}(x)\r]$$
and
$$\R=\frac{(\a+2)_{n-1}(\a+\b+2)_{n-1}}{(\b+1)_{n-1}\,n!}
\l[(n+\a+1)\P-(\a+1)P_n^{(\a+1,\b-1)}(x)\r].$$
By using these representations, (\ref{rel5}) and (\ref{diff})
we find that
\bea & &(1-x^2)D^2\Q+\l[\b-\a-(\a+\b+2)x\r]D\Q+n(n+\a+\b+1)\Q\nn
&=&-\frac{(\b+1)_n(\a+\b+2)_{n-1}}{(\a+1)_{n-1}\,n!}
\l[(1-x^2)D^2P_n^{(\a-1,\b+1)}(x)\r.\nn
& &{}\hspace{5cm}{}+\l[\b-\a-(\a+\b+2)x\r]DP_n^{(\a-1,\b+1)}(x)\nn
& &{}\hspace{7cm}\l.{}+n(n+\a+\b+1)P_n^{(\a-1,\b+1)}(x)\r]\nn
&=&\frac{(\b+1)_n(\a+\b+2)_{n-1}}{(\a+1)_{n-1}\,n!}2DP_n^{(\a-1,\b+1)}(x)\nn
&=&\frac{(\b+1)_n(\a+\b+2)_n}{(\a+1)_{n-1}\,n!}P_{n-1}^{(\a,\b+2)}(x),
\;n=1,2,3,\ldots\n\eea
and
\bea & &(1-x^2)D^2\R+\l[\b-\a-(\a+\b+2)x\r]D\R+n(n+\a+\b+1)\R\nn
&=&-\frac{(\a+1)_n(\a+\b+2)_{n-1}}{(\b+1)_{n-1}\,n!}
\l[(1-x^2)D^2P_n^{(\a+1,\b-1)}(x)\r.\nn
& &{}\hspace{5cm}{}+\l[\b-\a-(\a+\b+2)x\r]DP_n^{(\a+1,\b-1)}(x)\nn
& &{}\hspace{7cm}\l.{}+n(n+\a+\b+1)P_n^{(\a+1,\b-1)}(x)\r]\nn
&=&-\frac{(\a+1)_n(\a+\b+2)_{n-1}}{(\b+1)_{n-1}\,n!}2DP_n^{(\a+1,\b-1)}(x)\nn
&=&-\frac{(\a+1)_n(\a+\b+2)_n}{(\b+1)_{n-1}\,n!}P_{n-1}^{(\a+2,\b)}(x),
\;n=1,2,3,\ldots.\n\eea
Finally we have
$$D\l[(x^2-1)D^2\P\r]=(x^2-1)D^3\P+2xD^2\P,\;\ndots$$
and
\bea & &D^2\l[(x^2-1)D^2\P\r]\nn
& &{}\hspace{1cm}=(x^2-1)D^4\P+4xD^3\P+2D^2\P,\;\ndots,\n\eea
which implies by using (\ref{rel5}) that
\bea & &(1-x^2)D^2\l[(x^2-1)D^2\P\r]
+\l[\b-\a-(\a+\b+2)x\r]D\l[(x^2-1)D^2\P\r]\nn
& &{}\hspace{1cm}{}+n(n+\a+\b+1)(x^2-1)D^2\P\nn
&=&(x^2-1)\l[(1-x^2)D^4\P+\l[\b-\a-(\a+\b+6)x\r]D^3\P\r.\nn
& &{}\hspace{7cm}\l.{}+(n-2)(n+\a+\b+3)D^2\P\r]\nn
& &{}\hspace{1cm}{}+2\l[(\b+1)(x-1)-(\a+1)(x+1)\r]D^2\P\nn
&=&2\l[(\b+1)(x-1)-(\a+1)(x+1)\r]D^2\P,\;\ndots.\n\eea
Hence by using (\ref{defS1}) we find that
\bea & &(1-x^2)D^2\S+\l[\b-\a-(\a+\b+2)x\r]D\S+n(n+\a+\b+1)\S\nn
&=&\frac{2}{(\a+1)(\b+1)}\frac{(\a+\b+2)_n}{n!}
\frac{(\a+\b+2)_{n-1}}{(n-1)!}\nn
& &{}\hspace{1cm}\times\l[(\b+1)(x-1)-(\a+1)(x+1)\r]D^2\P,
\;n=1,2,3,\ldots.\n\eea

Since we demand that the differential equation (\ref{DV}) must hold for all
$M\ge 0$ and $N\ge 0$ we view the left-hand side of (\ref{invul}) as a
polynomial in $M$ and $N$ and conclude that all coefficients of this
polynomial must be equal to zero, hence we derive the following eight
systems of equations~:
$$\begin{array}{llll}
M~: & S_1=0 \hspace{3cm} & MN~: & S_5=0\\
M^2~: & S_2=0 & M^2N~: & S_6=0\\
N~: & S_3=0 & MN^2~: & S_7=0\\
N^2~: & S_4=0 & M^2N^2~: & S_8=0,
\end{array}$$
where $n=1,2,3,\ldots$ and
$$S_1=\sum_{i=0}^{\infty}a_i(x)D^i\P
+\frac{(\b+1)_n(\a+\b+2)_n}{(\a+1)_{n-1}\,n!}P_{n-1}^{(\a,\b+2)}(x),$$
$$S_2=\sum_{i=0}^{\infty}a_i(x)D^i\Q,$$
$$S_3=\sum_{i=0}^{\infty}b_i(x)D^i\P
-\frac{(\a+1)_n(\a+\b+2)_n}{(\b+1)_{n-1}\,n!}P_{n-1}^{(\a+2,\b)}(x),$$
$$S_4=\sum_{i=0}^{\infty}b_i(x)D^i\R,$$
\bea & &S_5=\sum_{i=0}^{\infty}a_i(x)D^i\R
+\sum_{i=0}^{\infty}b_i(x)D^i\Q+\sum_{i=0}^{\infty}c_i(x)D^i\P\nn
& &\hspace{1cm}{}+\frac{2}{(\a+1)(\b+1)}\frac{(\a+\b+2)_n}{n!}
\frac{(\a+\b+2)_{n-1}}{(n-1)!}\nn
& &{}\hspace{3cm}\times\l[(\b+1)(x-1)-(\a+1)(x+1)\r]D^2\P,\n\eea
$$S_6=\sum_{i=0}^{\infty}a_i(x)D^i\S
+\sum_{i=0}^{\infty}c_i(x)D^i\Q,$$
$$S_7=\sum_{i=0}^{\infty}b_i(x)D^i\S
+\sum_{i=0}^{\infty}c_i(x)D^i\R$$
and
$$S_8=\sum_{i=0}^{\infty}c_i(x)D^i\S.$$

By using (\ref{defQ2}) it follows from $S_1=0$ and $S_2=0$ that
$$\sum_{i=0}^{\infty}a_i(x)D^iP_{n-1}^{(\a,\b+1)}(x)=
\frac{(\b+2)_{n-1}(\a+\b+2)_n}{(\a+1)_{n-1}\,(n-1)!}P_{n-1}^{(\a,\b+2)}(x),
\;n=1,2,3,\ldots.$$
In view of (\ref{anul}) this is trivial for $n=1$. Hence, by shifting $n$
and using (\ref{anul}) and (\ref{rel3b}) we obtain
\bea & &\sum_{i=1}^{\infty}a_i(x)D^iP_n^{(\a,\b+1)}(x)\nn
&=&\frac{(\b+2)_n(\a+\b+2)_{n+1}}{(\a+1)_n\,n!}P_n^{(\a,\b+2)}(x)
-a_0(n+1,\a,\b)P_n^{(\a,\b+1)}(x)\nn
&=&\frac{(\b+2)_n(\a+\b+2)_{n+1}}{(\a+1)_n\,n!}P_n^{(\a,\b+2)}(x)
-\frac{(\b+3)_n(\a+\b+2)_{n+1}}{(\a+1)_n\,n!}P_n^{(\a,\b+1)}(x)\nn
&=&\frac{(\b+3)_{n-1}(\a+\b+2)_{n+1}}{(\a+1)_n\,n!}
\l[(\b+2)P_n^{(\a,\b+2)}(x)-(n+\b+2)P_n^{(\a,\b+1)}(x)\r]\nn
&=&-\frac{(\b+3)_{n-1}(\a+\b+2)_{n+1}}{(\a+1)_n\,n!}
(n+\a)\l(\frac{x+1}{2}\r)P_{n-1}^{(\a,\b+3)}(x)\nn
&=&-(\a+\b+2)\l(\frac{x+1}{2}\r)
\frac{(\b+3)_{n-1}(\a+\b+3)_n}{(\a+1)_{n-1}\,n!}P_{n-1}^{(\a,\b+3)}(x),
\;n=1,2,3,\ldots.\n\eea
Note that this system of equations has the form (\ref{inv1}). Hence by using
(\ref{opl1}) we conclude that
\bea & &a_i(\a,\b,x)=-(\a+\b+2)\l(\frac{x+1}{2}\r)2^i\nn
& &{}\hspace{3cm}{}\times\sum_{j=1}^i\frac{\a+\b+2j+2}{(\a+\b+j+2)_{i+1}}
\frac{(\b+3)_{j-1}(\a+\b+3)_j}{(\a+1)_{j-1}\,j!}\nn
& &{}\hspace{5cm}{}\times
P_{i-j}^{(-\a-i-1,-\b-i-2)}(x)P_{j-1}^{(\a,\b+3)}(x),\;i=1,2,3,\ldots.\n\eea
In the same way we obtain from $S_3=0$ and $S_4=0$ by using (\ref{defR2}),
(\ref{bnul}), (\ref{rel3a}), (\ref{inv1}) and (\ref{opl1})
\bea & &b_i(\a,\b,x)=-(\a+\b+2)\l(\frac{x-1}{2}\r)2^i\nn
& &{}\hspace{3cm}{}\times\sum_{j=1}^i\frac{\a+\b+2j+2}{(\a+\b+j+2)_{i+1}}
\frac{(\a+3)_{j-1}(\a+\b+3)_j}{(\b+1)_{j-1}\,j!}\nn
& &{}\hspace{5cm}{}\times
P_{i-j}^{(-\a-i-2,-\b-i-1)}(x)P_{j-1}^{(\a+3,\b)}(x),\;i=1,2,3,\ldots,\n\eea
but this is not really necessary in view of (\ref{symab}).

In order to prove (\ref{b}) we apply the definition (\ref{defJac2}) to
$P_{i-j-1}^{(-\a-i-2,-\b-i-1)}(x)$ and the definition (\ref{defJac1}) to
$P_j^{(\a+3,\b)}(x)$ to find by changing the order of summations and by
using the summation formula (\ref{sum2})
\bea & &\sum_{j=1}^i\frac{\a+\b+2j+2}{(\a+\b+j+2)_{i+1}}
\frac{(\a+3)_{j-1}(\a+\b+3)_j}{(\b+1)_{j-1}\,j!}
P_{i-j}^{(-\a-i-2,-\b-i-1)}(x)P_{j-1}^{(\a+3,\b)}(x)\nn
&=&\sum_{j=0}^{i-1}\frac{\a+\b+2j+4}{(\a+\b+j+3)_{i+1}}
\frac{(\a+3)_j(\a+\b+3)_{j+1}}{(\b+1)_j\,(j+1)!}
P_{i-j-1}^{(-\a-i-2,-\b-i-1)}(x)P_j^{(\a+3,\b)}(x)\nn
&=&\sum_{j=0}^{i-1}\frac{\a+\b+2j+4}{(\a+\b+j+3)_{i+1}}
\frac{(\a+3)_j(\a+\b+3)_{j+1}}{(\b+1)_j\,(j+1)!}\nn
& &{}\hspace{1cm}\times (-1)^{i-j-1}\sum_{k=0}^{i-j-1}
\frac{(\a+\b+i+j-k+4)_k}{k!}\frac{(\a+j+3)_{i-j-k-1}}{(i-j-k-1)!}
\l(\frac{x-1}{2}\r)^k\nn
& &{}\hspace{5cm}\times\sum_{m=0}^j\frac{(\a+\b+j+4)_m}{m!}
\frac{(\a+m+4)_{j-m}}{(j-m)!}\l(\frac{x-1}{2}\r)^m\nn
&=&\sum_{j=0}^{i-1}\sum_{k=0}^{i-j-1}\sum_{m=0}^j
(\a+\b+2j+4)(-1)^{i-j-1}\l(\frac{x-1}{2}\r)^{k+m}\nn
& &{}\hspace{1cm}\times\frac{(\a+3)_{i-k-1}(\a+\b+3)_{j+m+1}(\a+m+4)_{j-m}}
{(\a+\b+j+3)_{i-k+1}(\b+1)_j\,(j+1)!\,k!\,(i-j-k-1)!\,m!\,(j-m)!}\nn
&=&\sum_{k=0}^{i-1}\sum_{m=0}^{i-k-1}\sum_{j=0}^{i-k-m-1}
(\a+\b+2j+2m+4)(-1)^{i-j-m-1}\l(\frac{x-1}{2}\r)^{k+m}\nn
& &{}\hspace{1cm}\times\frac{(\a+3)_{i-k-1}(\a+\b+3)_{j+2m+1}(\a+m+4)_j}
{(\a+\b+j+m+3)_{i-k+1}(\b+1)_{j+m}\,(j+m+1)!\,k!\,(i-j-k-m-1)!\,m!\,j!}\nn
&=&\sum_{k=0}^{i-1}\sum_{m=0}^{i-k-1}\frac{(\a+3)_{i-k-1}(\a+\b+3)_{2m+1}}
{(\a+\b+m+3)_{i-k+1}(\b+1)_m\,(m+1)!\,k!\,(i-k-m-1)!\,m!}\nn
& &{}\hspace{9cm}\times(-1)^{i-m-1}\l(\frac{x-1}{2}\r)^{k+m}\nn
& &{}\hspace{1cm}\times\sum_{j=0}^{i-k-m-1}
\frac{(-i+k+m+1)_j(\a+\b+2m+4)_j(\a+m+4)_j(\a+\b+m+3)_j}
{(\a+\b+i-k+m+4)_j(\b+m+1)_j(m+2)_j\,j!}\nn
& &{}\hspace{9cm}\times(\a+\b+2j+2m+4)\nn
&=&\sum_{k=0}^{i-1}\sum_{m=0}^{i-k-1}\frac{(\a+3)_{i-k-1}(\a+\b+3)_{2m+1}}
{(\a+\b+m+3)_{i-k+1}(\b+1)_m\,(m+1)!\,k!\,(i-k-m-1)!\,m!}\nn
& &{}\hspace{1cm}\times\frac{(\a+\b+2m+4)_{i-k-m}(-\a-2)_{i-k-m-1}}
{(\b+m+1)_{i-k-m-1}(m+2)_{i-k-m-1}}(-1)^{i-m-1}\l(\frac{x-1}{2}\r)^{k+m}\nn
&=&\sum_{k=0}^{i-1}\sum_{m=0}^{i-k-1}(-1)^{i-m-1}
\frac{(\a+3)_{i-k-1}(\a+\b+3)_m(-\a-2)_{i-k-m-1}}
{(\b+1)_{i-k-1}\,(i-k)!\,k!\,(i-k-m-1)!\,m!}\l(\frac{x-1}{2}\r)^{k+m}\nn
&=&\sum_{m=0}^{i-1}\sum_{\ell=m}^{i-1}(-1)^{i-m-1}
\frac{(\a+3)_{i-\ell-1+m}(\a+\b+3)_m(-\a-2)_{i-\ell-1}}
{(\b+1)_{i-\ell-1+m}\,(i-\ell+m)!\,(\ell-m)!\,(i-\ell-1)!\,m!}
\l(\frac{x-1}{2}\r)^{\ell}\nn
&=&(-1)^{i-1}\sum_{\ell=0}^{i-1}\frac{(\a+3)_{i-\ell-1}(-\a-2)_{i-\ell-1}}
{(\b+1)_{i-\ell-1}\,(i-\ell)!\,\ell!\,(i-\ell-1)!}
\l(\frac{x-1}{2}\r)^{\ell}\nn
& &{}\hspace{3cm}\times\sum_{m=0}^{\ell}
\frac{(-\ell)_m(\a+i-\ell+2)_m(\a+\b+3)_m}{(\b+i-\ell)_m(i-\ell+1)_m\,m!},
\;i=1,2,3,\ldots.\n\eea
Hence
\bea & &b_i(\a,\b,x)=(\a+\b+2)(-2)^i
\sum_{\ell=0}^{i-1}\frac{(\a+3)_{i-\ell-1}(-\a-2)_{i-\ell-1}}
{(\b+1)_{i-\ell-1}\,(i-\ell)!\,(i-\ell-1)!\,\ell!}\nn
& &{}\hspace{3cm}\times
\hyp{3}{2}{-\ell,\a+i-\ell+2,\a+\b+3}{\b+i-\ell,i-\ell+1}{1}
\l(\frac{x-1}{2}\r)^{\ell+1},\;i=1,2,3\ldots,\n\eea
which proves (\ref{b}). The proof of (\ref{a}) is similar, but it is easier
to use (\ref{symab}) since then (\ref{a}) follows easily from (\ref{b}).

The computation of the coefficients $\{c_i(x)\}_{i=1}^{\infty}$ is more
difficult. First we set $n=1$ into $S_6=0$ or $S_7=0$. Since we have
$S_1^{(\a,\b)}(x)=0$ from (\ref{defS1}) and $c_0(1,\a,\b)=0$ we conclude
that
$$c_1(x)DQ_1^{(\a,\b)}(x)=0\;\textrm{ and }\;c_1(x)DR_1^{(\a,\b)}(x)=0.$$
By using (\ref{defQ1}), (\ref{defR1}) and (\ref{diff}) we find that
$$Q_1^{(\a,\b)}(x)=(\a+\b+2)\l(\frac{x+1}{2}\r)\;\textrm{ and }\;
R_1^{(\a,\b)}(x)=(\a+\b+2)\l(\frac{x-1}{2}\r).$$
Hence
$$DQ_1^{(\a,\b)}(x)=DR_1^{(\a,\b)}(x)=\frac{\a+\b+2}{2}\ne 0,$$
which implies that $c_1(x)=c_1(\a,\b,x)=0$.

Now we consider the system of equations $S_8=0$. Since
$S_1^{(\a,\b)}(x)=0$ the case $n=1$ is trivial. Now we use (\ref{defS1})
and (\ref{diff}) to find that
$$\S=\frac{1}{(\a+1)(\b+1)}\frac{(\a+\b+2)_n(\a+\b+2)_{n+1}}{4\,n!\,(n-1)!}
(x^2-1)P_{n-2}^{(\a+2,\b+2)}(x),\;n=2,3,4,\ldots.$$
By using the fact that
$$\frac{1}{(\a+1)(\b+1)}
\frac{(\a+\b+2)_n(\a+\b+2)_{n+1}}{4\,n!\,(n-1)!}\ne 0,\;n=2,3,4,\ldots$$
we conclude that
$$\sum_{i=0}^{\infty}c_i(x)D^i\l[(x^2-1)P_{n-2}^{(\a+2,\b+2)}(x)\r]=0,
\;n=2,3,4,\ldots.$$
Now we use the fact that $c_1(x)=0$ to obtain by shifting $n$
$$\sum_{i=2}^{\infty}c_i(x)D^i\l[(x^2-1)P_n^{(\a+2,\b+2)}(x)\r]=
c_0(n+2,\a,\b)(1-x^2)P_n^{(\a+2,\b+2)}(x),\;\ndots.$$
Note that for $\ndots$ we have
\bea & &D^i\l[(x^2-1)\P\r]=(x^2-1)D^i\P+2ixD^{i-1}\P\nn
& &{}\hspace{7cm}{}+i(i-1)D^{i-2}\P,\;i=2,3,4,\ldots.\n\eea
Hence we have
\be\la{sys}\sum_{i=0}^{\infty}C_i(x)D^iP_n^{(\a+2,\b+2)}(x)=
c_0(n+2,\a,\b)(1-x^2)P_n^{(\a+2,\b+2)}(x),\;\ndots,\ee
where
$$C_i(x)=\l\{\ba{ll}2c_2(x), & i=0\\[5mm]
4xc_2(x)+6c_3(x), & i=1\\[5mm]
(x^2-1)c_i(x)+2(i+1)xc_{i+1}(x)+(i+1)(i+2)c_{i+2}(x), & i=2,3,4,\ldots.
\ea\r.$$
Note that the system of equations (\ref{sys}) has the form (\ref{inv2}). So
we may apply (\ref{opl2}) and use (\ref{cnul}) to conclude that for
$i=0,1,2,\ldots$ we have
\bea C_i(x)&=&(1-x^2)2^i\sum_{j=0}^i\frac{\a+\b+2j+5}{(\a+\b+j+5)_{i+1}}
c_0(j+2,\a,\b)P_{i-j}^{(-\a-i-3,-\b-i-3)}(x)P_j^{(\a+2,\b+2)}(x)\nn
&=&\la{C1}\frac{(\a+\b+2)^2(\a+\b+3)}{(\a+1)(\b+1)}(1-x^2)2^i\nn
& &{}\hspace{1cm}\times\sum_{j=0}^i\frac{\a+\b+2j+5}{(\a+\b+j+5)_{i+1}}
\frac{(\a+\b+4)_{j+1}(\a+\b+4)_j}{(j+1)!\,j!}\nn
& &{}\hspace{5cm}\times
P_{i-j}^{(-\a-i-3,-\b-i-3)}(x)P_j^{(\a+2,\b+2)}(x).\eea

As before we apply the definition (\ref{defJac2}) to
$P_{i-j}^{(-\a-i-3,-\b-i-3)}(x)$ and the definition (\ref{defJac1}) to
$P_j^{(\a+2,\b+2)}(x)$ to find by changing the order of summations and by
using the summation formula (\ref{sum2})
\bea & &\sum_{j=0}^i\frac{\a+\b+2j+5}{(\a+\b+j+5)_{i+1}}
\frac{(\a+\b+4)_{j+1}(\a+\b+4)_j}{(j+1)!\,j!}
P_{i-j}^{(-\a-i-3,-\b-i-3)}(x)P_j^{(\a+2,\b+2)}(x)\nn
&=&\sum_{j=0}^i\frac{\a+\b+2j+5}{(\a+\b+j+5)_{i+1}}
\frac{(\a+\b+4)_{j+1}(\a+\b+4)_j}{(j+1)!\,j!}\nn
& &{}\hspace{1cm}\times (-1)^{i-j}\sum_{k=0}^{i-j}
\frac{(\a+\b+i+j-k+6)_k}{k!}\frac{(\a+j+3)_{i-j-k}}{(i-j-k)!}
\l(\frac{x-1}{2}\r)^k\nn
& &{}\hspace{3cm}\times\sum_{m=0}^j\frac{(\a+\b+j+5)_m}{m!}
\frac{(\a+m+3)_{j-m}}{(j-m)!}\l(\frac{x-1}{2}\r)^m\nn
&=&\sum_{j=0}^i\sum_{k=0}^{i-j}\sum_{m=0}^j
(\a+\b+2j+5)(-1)^{i-j}\l(\frac{x-1}{2}\r)^{k+m}\nn
& &{}\hspace{1cm}\times\frac{(\a+\b+4)_{j+m+1}(\a+\b+4)_j(\a+m+3)_{i-k-m}}
{(\a+\b+j+5)_{i-k+1}\,(j+1)!\,j!\,k!\,(i-j-k)!\,m!\,(j-m)!}\nn
&=&\sum_{k=0}^i\sum_{m=0}^{i-k}\sum_{j=0}^{i-k-m}
(\a+\b+2j+2m+5)(-1)^{i-j-m}\l(\frac{x-1}{2}\r)^{k+m}\nn
& &{}\hspace{1cm}\times
\frac{(\a+\b+4)_{j+2m+1}(\a+\b+4)_{j+m}(\a+m+3)_{i-k-m}}
{(\a+\b+j+m+5)_{i-k+1}\,(j+m+1)!\,(j+m)!\,k!\,(i-j-k-m)!\,m!\,j!}\nn
&=&\sum_{k=0}^i\sum_{m=0}^{i-k}
\frac{(\a+\b+4)_{2m+1}(\a+\b+4)_m(\a+m+3)_{i-k-m}}
{(\a+\b+m+5)_{i-k+1}\,(m+1)!\,m!\,k!\,(i-k-m)!\,m!}\nn
& &{}\hspace{9cm}\times(-1)^{i-m}\l(\frac{x-1}{2}\r)^{k+m}\nn
& &{}\hspace{1cm}\times\sum_{j=0}^{i-k-m}
\frac{(-i+k+m)_j(\a+\b+2m+5)_j(\a+\b+m+4)_j(\a+\b+m+5)_j}
{(\a+\b+i-k+m+6)_j(m+2)_j(m+1)_j\,j!}\nn
& &{}\hspace{9cm}\times(\a+\b+2j+2m+5)\nn
&=&\sum_{k=0}^i\sum_{m=0}^{i-k}
\frac{(\a+\b+4)_{2m+1}(\a+\b+4)_m(\a+m+3)_{i-k-m}}
{(\a+\b+m+5)_{i-k+1}\,(m+1)!\,m!\,k!\,(i-k-m)!\,m!}\nn
& &{}\hspace{1cm}\times\frac{(\a+\b+2m+5)_{i-k-m+1}(-\a-\b-3)_{i-k-m}}
{(m+2)_{i-k-m}(m+1)_{i-k-m}}(-1)^{i-m}\l(\frac{x-1}{2}\r)^{k+m}\nn
&=&\sum_{k=0}^i\sum_{m=0}^{i-k}
\frac{(\a+\b+4)_{m+1}(\a+\b+4)_m(\a+m+3)_{i-k-m}(-\a-\b-3)_{i-k-m}}
{(i-k+1)!\,(i-k)!\,k!\,(i-k-m)!\,m!}\nn
& &{}\hspace{9cm}\times (-1)^{i-m}\l(\frac{x-1}{2}\r)^{k+m}\nn
&=&\sum_{m=0}^i\sum_{\ell=m}^i
\frac{(\a+\b+4)_{m+1}(\a+\b+4)_m(\a+m+3)_{i-\ell}(-\a-\b-3)_{i-\ell}}
{(i-\ell+m+1)!\,(i-\ell+m)!\,(\ell-m)!\,(i-\ell)!\,m!}\nn
& &{}\hspace{9cm}\times (-1)^{i-m}\l(\frac{x-1}{2}\r)^{\ell}\nn
&=&(\a+\b+4)(-1)^i\sum_{\ell=0}^i\frac{(\a+3)_{i-\ell}(-\a-\b-3)_{i-\ell}}
{(i-\ell+1)!\,(i-\ell)!\,\ell!\,(i-\ell)!}
\l(\frac{x-1}{2}\r)^{\ell}\nn
& &{}\hspace{1cm}\times\sum_{m=0}^{\ell}
\frac{(-\ell)_m(\a+\b+5)_m(\a+\b+4)_m(\a+i-\ell+3)_m}
{(\a+3)_m(i-\ell+2)_m(i-\ell+1)_m\,m!},\;i=0,1,2,\ldots.\n\eea
Hence for $i=0,1,2,\ldots$ we have
\bea\la{C2}& &C_i(x)=\frac{(\a+\b+2)^2(\a+\b+3)(\a+\b+4)}{(\a+1)(\b+1)}
(1-x^2)(-2)^i\nn
& &{}\hspace{2cm}\times\sum_{\ell=0}^i
\frac{(\a+3)_{i-\ell}(-\a-\b-3)_{i-\ell}}
{(i-\ell+1)!\,(i-\ell)!\,\ell!\,(i-\ell)!}\nn
& &{}\hspace{3cm}\times
\hyp{4}{3}{-\ell,\a+\b+5,\a+\b+4,\a+i-\ell+3}{\a+3,i-\ell+2,i-\ell+1}{1}
\l(\frac{x-1}{2}\r)^{\ell}.\eea

Now we have by using (\ref{A2}) and (\ref{B2})
\bea c_i(x)&=&\frac{1}{2\,i!}\sum_{j=0}^{i-2}(-1)^{i-j}j!\,
\l[(x+1)^{i-j-1}-(x-1)^{i-j-1}\r]C_j(x)\nn
&=&\frac{1}{2\,i!}\sum_{j=0}^{i-2}(-1)^j(i-j-2)!\,
\l[(x+1)^{j+1}-(x-1)^{j+1}\r]C_{i-j-2}(x)\nn
&=&c_i^{(1)}(x)+c_i^{(2)}(x),\;i=2,3,4,\ldots,\n\eea
where
$$c_i^{(1)}(x)=c_i^{(1)}(\a,\b,x)=\frac{1}{2\,i!}\sum_{j=0}^{i-2}
(-1)^j(i-j-2)!\,(x+1)^{j+1}C_{i-j-2}(x),\;i=2,3,4,\ldots$$
and
$$c_i^{(2)}(x)=c_i^{(2)}(\a,\b,x)=\frac{1}{2\,i!}\sum_{j=0}^{i-2}
(-1)^{j+1}(i-j-2)!\,(x-1)^{j+1}C_{i-j-2}(x),\;i=2,3,4,\ldots.$$

Now we will prove (\ref{c1}) and (\ref{c2}). To do this we will
first prove (\ref{symc}), which is an easy consequence of the
symmetry formula (\ref{symJac}). If we write
$C_i(x)=C_i(\a,\b,x)$ this symmetry formula gives us
$$C_i(\a,\b,x)=(-1)^iC_i(\b,\a,-x),\;i=0,1,2,\ldots$$
in view of (\ref{C1}). Hence
\bea c_i^{(1)}(\a,\b,x)&=&\frac{1}{2\,i!}\sum_{j=0}^{i-2}
(-1)^j(i-j-2)!\,(x+1)^{j+1}C_{i-j-2}(\a,\b,x)\nn
&=&\frac{1}{2\,i!}\sum_{j=0}^{i-2}
(-1)^{i+j+1}(i-j-2)!\,(-x-1)^{j+1}C_{i-j-2}(\b,\a,-x)\nn
&=&(-1)^ic_i^{(2)}(\b,\a,-x),\;i=2,3,4,\ldots,\n\eea
which proves (\ref{symc}). In order to prove (\ref{c2}) we use
(\ref{C2}) and change the order of summations to find for
$i=2,3,4,\ldots$
\bea c_i^{(2)}(x)&=&\frac{1}{2\,i!}\sum_{j=0}^{i-2}
(-1)^{j+1}(i-j-2)!\,(x-1)^{j+1}C_{i-j-2}(x)\nn
&=&\frac{(\a+\b+2)^2(\a+\b+3)(\a+\b+4)}{(\a+1)(\b+1)}
(x^2-1)\frac{(-2)^{i-2}}{i!}\nn
& &{}\hspace{1cm}\times\sum_{j=0}^{i-2}\sum_{k=0}^{i-j-2}
\frac{(\a+3)_{i-j-k-2}(-\a-\b-3)_{i-j-k-2}(i-j-2)!}
{(i-j-k-1)!\,(i-j-k-2)!\,k!\,(i-j-k-2)!}\nn
& &{}\hspace{1cm}\times
\hyp{4}{3}{-k,\a+\b+5,\a+\b+4,\a+i-j-k+1}{\a+3,i-j-k,i-j-k-1}{1}
\l(\frac{x-1}{2}\r)^{j+k+1}\nn
&=&\frac{(\a+\b+2)^2(\a+\b+3)(\a+\b+4)}{(\a+1)(\b+1)}
(x^2-1)\frac{(-2)^{i-2}}{i!}\nn
& &{}\hspace{1cm}\times\sum_{j=0}^{i-2}\sum_{\ell=j}^{i-2}
\frac{(\a+3)_{i-\ell-2}(-\a-\b-3)_{i-\ell-2}(i-j-2)!}
{(i-\ell-1)!\,(i-\ell-2)!\,(\ell-j)!\,(i-\ell-2)!}\nn
& &{}\hspace{1cm}\times
\hyp{4}{3}{-\ell+j,\a+\b+5,\a+\b+4,\a+i-\ell+1}{\a+3,i-\ell,i-\ell-1}{1}
\l(\frac{x-1}{2}\r)^{\ell+1}\nn
&=&\frac{(\a+\b+2)^2(\a+\b+3)(\a+\b+4)}{(\a+1)(\b+1)}
(x^2-1)\frac{(-2)^{i-2}}{i!}\nn
& &{}\hspace{1cm}\times\sum_{\ell=0}^{i-2}\sum_{j=0}^{\ell}
\frac{(\a+3)_{i-\ell-2}(-\a-\b-3)_{i-\ell-2}(i-j-2)!}
{(i-\ell-1)!\,(i-\ell-2)!\,(\ell-j)!\,(i-\ell-2)!}\nn
& &{}\hspace{1cm}\times\sum_{m=0}^{\ell-j}
\frac{(-\ell+j)_m(\a+\b+5)_m(\a+\b+4)_m(\a+i-\ell+1)_m}
{(\a+3)_m(i-\ell)_m(i-\ell-1)_m\,m!}\l(\frac{x-1}{2}\r)^{\ell+1}\nn
&=&\frac{(\a+\b+2)^2(\a+\b+3)(\a+\b+4)}{(\a+1)(\b+1)}
(x^2-1)\frac{(-2)^{i-2}}{i!}\nn
& &{}\hspace{1cm}\times\sum_{\ell=0}^{i-2}\sum_{m=0}^{\ell}
\frac{(\a+3)_{i-\ell-2}(-\a-\b-3)_{i-\ell-2}}
{(i-\ell-1)!\,(i-\ell-2)!\,(i-\ell-2)!}\nn
& &{}\hspace{1cm}\times\frac{(\a+\b+5)_m(\a+\b+4)_m(\a+i-\ell+1)_m}
{(\a+3)_m(i-\ell)_m(i-\ell-1)_m\,m!}\l(\frac{x-1}{2}\r)^{\ell+1}\nn
& &{}\hspace{7cm}\times
\sum_{j=0}^{\ell-m}\frac{(i-j-2)!\,(-\ell+j)_m}{(\ell-j)!}.\n\eea
Now we use the Vandermonde summation formula (\ref{Van}) to obtain
\bea & &\sum_{j=0}^{\ell-m}\frac{(i-j-2)!\,(-\ell+j)_m}{(\ell-j)!}
=\frac{(i-2)!}{\ell!}(-\ell)_m\sum_{j=0}^{\ell-m}
\frac{(-\ell+m)_j}{(-i+2)_j}\nn
& &{}\hspace{1cm}=\frac{(i-2)!}{\ell!}(-\ell)_m\,\,
\hyp{2}{1}{-\ell+m,1}{-i+2}{1}=\frac{(i-2)!}{\ell!}(-\ell)_m
\frac{(-i+1)_{\ell-m}}{(-i+2)_{\ell-m}}\nn
& &{}\hspace{1cm}=\frac{(i-2)!}{\ell!}(-\ell)_m
\frac{(i-1)!\,(i-\ell+m-2)!}{(i-2)!\,(i-\ell+m-1)!}
=\frac{(i-1)!\,(i-\ell-2)!}{\ell!\,(i-\ell-1)!}
\frac{(-\ell)_m(i-\ell-1)_m}{(i-\ell)_m}.\n\eea
Hence
\bea c_i^{(2)}(x)&=&
\frac{(\a+\b+2)^2(\a+\b+3)(\a+\b+4)}{(\a+1)(\b+1)i}(x^2-1)(-2)^{i-2}\nn
& &{}\hspace{1cm}\times\sum_{\ell=0}^{i-2}\sum_{m=0}^{\ell}
\frac{(\a+3)_{i-\ell-2}(-\a-\b-3)_{i-\ell-2}}
{(i-\ell-1)!\,(i-\ell-2)!\,\ell!\,(i-\ell-1)!}\nn
& &{}\hspace{2cm}\times\frac{(-\ell)_m(\a+\b+5)_m(\a+\b+4)_m(\a+i-\ell+1)_m}
{(\a+3)_m(i-\ell)_m(i-\ell)_m\,m!}\l(\frac{x-1}{2}\r)^{\ell+1}\nn
&=&\frac{(\a+\b+2)^2(\a+\b+3)(\a+\b+4)}{(\a+1)(\b+1)i}(x^2-1)(-2)^{i-2}\nn
& &{}\hspace{1cm}\times\sum_{\ell=0}^{i-2}
\frac{(\a+3)_{i-\ell-2}(-\a-\b-3)_{i-\ell-2}}
{(i-\ell-1)!\,(i-\ell-2)!\,\ell!\,(i-\ell-1)!}\nn
& &{}\hspace{2cm}\times\hyp{4}{3}{-\ell,\a+\b+5,\a+\b+4,\a+i-\ell+1}
{\a+3,i-\ell,i-\ell}{1}\l(\frac{x-1}{2}\r)^{\ell+1}\n\eea
for $i=2,3,4,\ldots$, which proves (\ref{c2}).

Hence we have proved (\ref{c}), (\ref{c1}) and (\ref{c2}).

\section{The order of the differential equation}

For $\a>-1$, $\b>-1$, $M\ge 0$ and $N\ge 0$ the generalized Jacobi
polynomials $\set{\GP}$ satisfy a unique differential equation of the form
(\ref{DV}), where the coefficients are given by (\ref{anul}), (\ref{bnul}),
(\ref{cnul}), (\ref{a}), (\ref{b}), (\ref{c}), (\ref{c1}) and (\ref{c2}).

First of all we remark that
$$a_i(\a,\b,x)=\sum_{j=1}^ik_{i,j}^{(a)}(\a,\b)(x+1)^j,\;i=1,2,3,\ldots,$$
where
$$k_{i,1}^{(a)}(\a,\b)=-(\a+\b+2)2^{i-1}
\frac{(\b+3)_{i-1}(-\b-2)_{i-1}}{(\a+1)_{i-1}\,i!\,(i-1)!},
\;i=1,2,3,\ldots.$$
Since $\a>-1$ and $\b>-1$ we conclude that $k_{i,1}^{(a)}(\a,\b)$ only vanishes
if $\b\in\{0,1,2,\ldots\}$ and $i\ge\b+4$.
In the same way we have
$$b_i(\a,\b,x)=\sum_{j=1}^ik_{i,j}^{(b)}(\a,\b)(x-1)^j,\;i=1,2,3,\ldots,$$
where
$$k_{i,1}^{(b)}(\a,\b)=-(\a+\b+2)(-2)^{i-1}
\frac{(\a+3)_{i-1}(-\a-2)_{i-1}}{(\b+1)_{i-1}\,i!\,(i-1)!},
\;i=1,2,3,\ldots.$$
Hence, $k_{i,1}^{(b)}(\a,\b)$ only vanishes if $\a\in\{0,1,2,\ldots\}$ and
$i\ge\a+4$.

Now we will prove (\ref{afbrb}). So let $\a\in\{0,1,2,\ldots\}$ and consider
$b_i(\a,\b,x)$ given by (\ref{b}). Suppose that $i\ge 2\a+4$. Then we have
$$(-\a-2)_{i-\ell-1}=0\;\textrm{ for }\;i-\ell\ge\a+4.$$
Suppose that $i-\ell\le\a+3$, then we have $\ell\ge i-\a-3\ge\a+1$. Hence by
using (\ref{trans2}) and the Vandermonde summation formula (\ref{Van}) we
find that
\bea & &\hyp{3}{2}{-\ell,\a+\b+3,\a+i-\ell+2}{\b+i-\ell,i-\ell+1}{1}\nn
&=&\sum_{n=0}^{\ell}(-1)^n\frac{(-\ell)_n(\a+\b+3)_n(-\a-1)_n}
{(\b+i-\ell)_n(i-\ell+1)_n\,n!}\,
\hyp{2}{1}{n-\ell,n+\a+\b+3}{n+\b+i-\ell}{1}\nn
&=&\frac{1}{(\b+i-\ell)_{\ell}}\sum_{n=0}^{\a+1}(-1)^n
\frac{(-\ell)_n(\a+\b+3)_n(-\a-1)_n(i-\ell-\a-3)_{\ell-n}}
{(i-\ell+1)_n\,n!}.\n\eea
Since $i-\ell\le\a+3$ we have $i-\ell-\a-3\le 0$. Hence
$(i-\ell-\a-3)_{\ell-n}=0$ for $\ell-n\ge-i+\ell+\a+4$ or $i\ge n+\a+4$.
This implies that $b_i(\a,\b,x)=0$ if $i>n+\a+3$ for all
$n\in\{0,1,2,\ldots,\a+1\}$, hence for $i>\a+1+\a+3=2\a+4$.

In the same way we obtain (\ref{afbra}).

Now we will prove (\ref{akop}) and (\ref{bkop}). Suppose that
$\a\in\{0,1,2,\ldots\}$ and $i=2\a+4$. Then we have
$$(i-\ell-\a-3)_{\ell-n}=(\a+1-\ell)_{\ell-n}=0\;\textrm{ for }\;n\le\a.$$
Hence
\bea & &b_{2\a+4}(\a,\b,x)\nn
&=&(\a+\b+2)(-2)^{2\a+4}\sum_{\ell=\a+1}^{2\a+3}
\frac{(\a+3)_{2\a+3-\ell}(-\a-2)_{2\a+3-\ell}}
{(\b+1)_{2\a+3-\ell}(2\a+4-\ell)!\,(2\a+3-\ell)!\,\ell!}
\l(\frac{x-1}{2}\r)^{\ell+1}\nn
& &{}\hspace{1cm}\times\frac{(-1)^{\a+1}}{(2\a+\b+4-\ell)_{\ell}}
\frac{(-\ell)_{\a+1}(\a+\b+3)_{\a+1}(-\a-1)_{\a+1}(\a+1-\ell)_{\ell-\a-1}}
{(2\a+5-\ell)_{\a+1}(\a+1)!}\nn
&=&(\a+\b+2)(-2)^{2\a+4}\frac{(\a+\b+3)_{\a+1}}{(\b+1)_{2\a+3}(\a+2)!}
\sum_{\ell=\a+1}^{2\a+3}(-1)^{\ell}
\frac{(-\a-2)_{2\a+3-\ell}}{(2\a+3-\ell)!}\l(\frac{x-1}{2}\r)^{\ell+1}\nn
&=&-\frac{2^{2\a+4}}{(\b+1)_{\a+1}(\a+2)!}\sum_{\ell=0}^{\a+2}
(-1)^{\ell}\frac{(-\a-2)_{\ell}}{\ell!}\l(\frac{x-1}{2}\r)^{2\a+4-\ell}\nn
&=&-\frac{2^{2\a+4}}{(\b+1)_{\a+1}(\a+2)!}
\l(\frac{x-1}{2}\r)^{\a+2}\l(\frac{x+1}{2}\r)^{\a+2}
=-\frac{(x^2-1)^{\a+2}}{(\b+1)_{\a+1}(\a+2)!},\n\eea
which proves (\ref{bkop}). The proof of (\ref{akop}) is similar.

In order to prove (\ref{afbrc}) we first consider $c_i^{(2)}(\a,\b,x)$
given by (\ref{c2}). Let $\a,\b\in\{0,1,2,\ldots\}$ and suppose that
$i\ge 2\a+2\b+7$. Then we have
$$(-\a-\b-3)_{i-\ell-2}=0\;\textrm{ for }\;i-\ell\ge\a+\b+6.$$
Suppose that $i-\ell\le\a+\b+5$, then we have $\ell\ge i-\a-\b-5\ge\a+\b+2$.
Hence by using (\ref{trans1}), (\ref{trans2}) and the Vandermonde summation
formula (\ref{Van}) we find that
\bea & &
\hyp{4}{3}{-\ell,\a+\b+5,\a+\b+4,\a+i-\ell+1}{\a+3,i-\ell,i-\ell}{1}\nn
&=&\sum_{n=0}^{\ell}(-1)^n\frac{(-\ell)_n(\a+\b+5)_n(\a+\b+4)_n(-\a-1)_n}
{(\a+3)_n(i-\ell)_n(i-\ell)_n\,n!}\nn
& &{}\hspace{3cm}\times
\hyp{3}{2}{n-\ell,n+\a+\b+5,n+\a+\b+4}{n+i-\ell,n+\a+3}{1}\nn
&=&\sum_{n=0}^{\a+1}(-1)^n\frac{(-\ell)_n(\a+\b+5)_n(\a+\b+4)_n(-\a-1)_n}
{(\a+3)_n(i-\ell)_n(i-\ell)_n\,n!}\nn
& &{}\hspace{3cm}\times\sum_{k=0}^{\ell-n}(-1)^k
\frac{(n-\ell)_k(n+\a+\b+5)_k(-\b-1)_k}{(n+i-\ell)_k(n+\a+3)_k\,k!}\nn
& &{}\hspace{7cm}\times\hyp{2}{1}{n+k-\ell,n+k+\a+\b+5}{n+k+i-\ell}{1}\nn
&=&\sum_{n=0}^{\a+1}(-1)^n\frac{(-\ell)_n(\a+\b+5)_n(\a+\b+4)_n(-\a-1)_n}
{(\a+3)_n(i-\ell)_n(i-\ell)_n\,n!}\nn
& &{}\hspace{3cm}\times\sum_{k=0}^{\b+1}(-1)^k
\frac{(n-\ell)_k(n+\a+\b+5)_k(-\b-1)_k}{(n+i-\ell)_k(n+\a+3)_k\,k!}\nn
& &{}\hspace{7cm}\times
\frac{(i-\ell-\a-\b-5)_{\ell-n-k}}{(n+k+i-\ell)_{\ell-n-k}}\nn
&=&\frac{1}{(i-\ell)_{\ell}}\sum_{n=0}^{\a+1}\sum_{k=0}^{\b+1}(-1)^{n+k}
\frac{(-\ell)_{n+k}(\a+\b+5)_{n+k}(\a+\b+4)_n}
{(\a+3)_{n+k}(i-\ell)_n\,n!\,k!}\nn
& &{}\hspace{5cm}\times(-\a-1)_n(-\b-1)_k(i-\ell-\a-\b-5)_{\ell-n-k}.\n\eea
Since $i-\ell\le\a+\b+5$ we have $i-\ell-\a-\b-5\le 0$. Hence
$(i-\ell-\a-\b-5)_{\ell-n-k}=0$ for $\ell-n-k\ge-i+\ell+\a+\b+6$ or
$i\ge n+k+\a+\b+6$. This implies that $c_i^{(2)}(\a,\b,x)=0$ if $i>n+k+\a+\b+5$
for all $n\in\{0,1,2,\ldots,\a+1\}$ and
$k\in\{0,1,2,\ldots,\b+1\}$, hence for $i>\a+1+\b+1+\a+\b+5=2\a+2\b+7$.

In the same way we find that $c_i^{(1)}(\a,\b,x)=0$ for $i>2\a+2\b+7$.

Suppose that $\a,\b\in\{0,1,2,\ldots\}$ and $i=2\a+2\b+7$. Then we have
$$(i-\ell-\a-\b-5)_{\ell-n-k}=(\a+\b+2-\ell)_{\ell-n-k}=0\;\textrm{ for }
\;n+k\le\a+\b+1.$$
Hence
\bea & &c_{2\a+2\b+7}^{(2)}(\a,\b,x)\nn
&=&\frac{(\a+\b+2)^2(\a+\b+3)(\a+\b+4)}{(\a+1)(\b+1)(2\a+2\b+7)}
(x^2-1)(-2)^{2\a+2\b+5}\nn
& &{}\times\sum_{\ell=\a+\b+2}^{2\a+2\b+5}
\frac{(\a+3)_{2\a+2\b+5-\ell}(-\a-\b-3)_{2\a+2\b+5-\ell}}
{(2\a+2\b+6-\ell)!\,(2\a+2\b+5-\ell)!\,(2\a+2\b+6-\ell)!\,\ell!}
\l(\frac{x-1}{2}\r)^{\ell+1}\nn
& &{}\hspace{1cm}\times\frac{(-1)^{\a+\b+2}}{(2\a+2\b+7-\ell)_{\ell}}
\frac{(-\ell)_{\a+\b+2}(\a+\b+5)_{\a+\b+2}(\a+\b+4)_{\a+1}}
{(\a+3)_{\a+\b+2}(2\a+2\b+7-\ell)_{\a+1}\,(\a+1)!\,(\b+1)!}\nn
& &{}\hspace{5cm}\times
(-\a-1)_{\a+1}(-\b-1)_{\b+1}(\a+\b+2-\ell)_{\ell-\a-\b-2}\nn
&=&\frac{(\a+\b+2)^2(\a+\b+3)(\a+\b+4)}{(\a+1)(\b+1)(2\a+2\b+7)}
(x^2-1)(-2)^{2\a+2\b+5}\nn
& &{}\hspace{1cm}\times\frac{(\a+\b+5)_{\a+\b+2}(\a+\b+4)_{\a+1}}
{(\a+3)_{\a+\b+2}(2\a+2\b+6)!\,(\a+2)!}\nn
& &{}\hspace{5cm}\times\sum_{\ell=\a+\b+2}^{2\a+2\b+5}(-1)^{\ell}
\frac{(-\a-\b-3)_{2\a+2\b+5-\ell}}{(2\a+2\b+5-\ell)!}
\l(\frac{x-1}{2}\r)^{\ell+1}\nn
&=&\frac{\a+\b+2}{(\a+1)(\b+1)(2\a+2\b+7)}
\frac{2^{2\a+2\b+5}}{(\a+\b+1)!\,(\a+\b+3)!}(x^2-1)\nn
& &{}\hspace{5cm}\times\sum_{\ell=0}^{\a+\b+3}(-1)^{\ell}
\frac{(-\a-\b-3)_{\ell}}{\ell!}\l(\frac{x-1}{2}\r)^{2\a+2\b+6-\ell}\nn
&=&\frac{\a+\b+2}{(\a+1)(\b+1)(2\a+2\b+7)}
\frac{2^{2\a+2\b+5}}{(\a+\b+1)!\,(\a+\b+3)!}(x^2-1)\nn
& &{}\hspace{5cm}\times
\l(\frac{x-1}{2}\r)^{\a+\b+3}\l(\frac{x+1}{2}\r)^{\a+\b+3}\nn
&=&\frac{\a+\b+2}{2(\a+1)(\b+1)(2\a+2\b+7)}
\frac{(x^2-1)^{\a+\b+4}}{(\a+\b+1)!\,(\a+\b+3)!}.\n\eea
Hence, because of the symmetry relation (\ref{symc}) we find that
$$c_{2\a+2\b+7}^{(1)}(\a,\b,x)=-c_{2\a+2\b+7}^{(2)}(\b,\a,-x)
=-c_{2\a+2\b+7}^{(2)}(\a,\b,x),$$
which implies that
$$c_{2\a+2\b+7}(\a,\b,x)=c_{2\a+2\b+7}^{(1)}(\a,\b,x)
+c_{2\a+2\b+7}^{(2)}(\a,\b,x)=0.$$
Hence we have proved (\ref{afbrc}).

In order to prove (\ref{ckop}) we first consider $C_i(\a,\b,x)=C_i(x)$
given by (\ref{C2}). We assume that $\a,\b\in\{0,1,2,\ldots\}$ and
$i\ge 2\a+2\b+6$. Then we have as before
$$(-\a-\b-3)_{i-\ell}=0\;\textrm{ for }\;i-\ell\ge\a+\b+4.$$
Suppose that $i-\ell\le\a+\b+3$, then we have $\ell\ge i-\a-\b-3\ge\a+\b+3$.
Hence by using (\ref{trans1}), (\ref{trans2}) and the Vandermonde summation
formula (\ref{Van}) we find that
\bea & &
\hyp{4}{3}{-\ell,\a+\b+5,\a+\b+4,\a+i-\ell+3}{\a+3,i-\ell+1,i-\ell+2}{1}\nn
&=&\sum_{n=0}^{\ell}(-1)^n\frac{(-\ell)_n(\a+\b+5)_n(\a+\b+4)_n(-\a-1)_n}
{(\a+3)_n(i-\ell+1)_n(i-\ell+2)_n\,n!}\nn
& &{}\hspace{3cm}\times
\hyp{3}{2}{n-\ell,n+\a+\b+5,n+\a+\b+4}{n+i-\ell+1,n+\a+3}{1}\nn
&=&\sum_{n=0}^{\a+1}(-1)^n\frac{(-\ell)_n(\a+\b+5)_n(\a+\b+4)_n(-\a-1)_n}
{(\a+3)_n(i-\ell+1)_n(i-\ell+2)_n\,n!}\nn
& &{}\hspace{3cm}\times\sum_{k=0}^{\ell-n}(-1)^k
\frac{(n-\ell)_k(n+\a+\b+5)_k(-\b-1)_k}{(n+i-\ell+1)_k(n+\a+3)_k\,k!}\nn
& &{}\hspace{7cm}\times\hyp{2}{1}{n+k-\ell,n+k+\a+\b+5}{n+k+i-\ell+1}{1}\nn
&=&\sum_{n=0}^{\a+1}(-1)^n\frac{(-\ell)_n(\a+\b+5)_n(\a+\b+4)_n(-\a-1)_n}
{(\a+3)_n(i-\ell+1)_n(i-\ell+2)_n\,n!}\nn
& &{}\hspace{3cm}\times\sum_{k=0}^{\b+1}(-1)^k
\frac{(n-\ell)_k(n+\a+\b+5)_k(-\b-1)_k}{(n+i-\ell+1)_k(n+\a+3)_k\,k!}\nn
& &{}\hspace{7cm}\times
\frac{(i-\ell-\a-\b-4)_{\ell-n-k}}{(n+k+i-\ell+1)_{\ell-n-k}}\nn
&=&\frac{1}{(i-\ell+1)_{\ell}}\sum_{n=0}^{\a+1}\sum_{k=0}^{\b+1}(-1)^{n+k}
\frac{(-\ell)_{n+k}(\a+\b+5)_{n+k}(\a+\b+4)_n}
{(\a+3)_{n+k}(i-\ell+2)_n\,n!\,k!}\nn
& &{}\hspace{5cm}\times(-\a-1)_n(-\b-1)_k(i-\ell-\a-\b-4)_{\ell-n-k}.\n\eea
Since $i-\ell\le\a+\b+3$ we have $i-\ell-\a-\b-4\le -1$. Hence
$(i-\ell-\a-\b-4)_{\ell-n-k}=0$ for $\ell-n-k\ge-i+\ell+\a+\b+5$ or
$i\ge n+k+\a+\b+5$. This implies that $C_i(x)=0$ if $i>n+k+\a+\b+4$
for all $n\in\{0,1,2,\ldots,\a+1\}$ and
$k\in\{0,1,2,\ldots,\b+1\}$, hence for $i>\a+1+\b+1+\a+\b+4=2\a+2\b+6$.

Suppose that $\a,\b\in\{0,1,2,\ldots\}$ and $i=2\a+2\b+6$. Then we have
$$(i-\ell-\a-\b-4)_{\ell-n-k}=(\a+\b+2-\ell)_{\ell-n-k}=0\;\textrm{ for }
\;n+k\le\a+\b+1.$$
Hence
\bea & &C_{2\a+2\b+6}(x)\nn
&=&\frac{(\a+\b+2)^2(\a+\b+3)(\a+\b+4)}{(\a+1)(\b+1)}
(1-x^2)(-2)^{2\a+2\b+6}\nn
& &{}\times\sum_{\ell=\a+\b+3}^{2\a+2\b+6}
\frac{(\a+3)_{2\a+2\b+6-\ell}(-\a-\b-3)_{2\a+2\b+6-\ell}}
{(2\a+2\b+7-\ell)!\,(2\a+2\b+6-\ell)!\,(2\a+2\b+6-\ell)!\,\ell!}
\l(\frac{x-1}{2}\r)^{\ell}\nn
& &{}\hspace{1cm}\times\frac{(-1)^{\a+\b+2}}{(2\a+2\b+7-\ell)_{\ell}}
\frac{(-\ell)_{\a+\b+2}(\a+\b+5)_{\a+\b+2}(\a+\b+4)_{\a+1}}
{(\a+3)_{\a+\b+2}(2\a+2\b+8-\ell)_{\a+1}\,(\a+1)!\,(\b+1)!}\nn
& &{}\hspace{5cm}\times
(-\a-1)_{\a+1}(-\b-1)_{\b+1}(\a+\b+2-\ell)_{\ell-\a-\b-2}\nn
&=&\frac{(\a+\b+2)^2(\a+\b+3)(\a+\b+4)}{(\a+1)(\b+1)}
(1-x^2)(-2)^{2\a+2\b+6}\nn
& &{}\hspace{1cm}\times\frac{(\a+\b+5)_{\a+\b+2}(\a+\b+4)_{\a+1}}
{(\a+3)_{\a+\b+2}(2\a+2\b+6)!\,(\a+2)!}\nn
& &{}\hspace{5cm}\times\sum_{\ell=\a+\b+3}^{2\a+2\b+6}(-1)^{\ell}
\frac{(-\a-\b-3)_{2\a+2\b+6-\ell}}{(2\a+2\b+6-\ell)!}
\l(\frac{x-1}{2}\r)^{\ell}\nn
&=&\frac{\a+\b+2}{(\a+1)(\b+1)}
\frac{2^{2\a+2\b+6}}{(\a+\b+1)!\,(\a+\b+3)!}(1-x^2)\nn
& &{}\hspace{5cm}\times\sum_{\ell=0}^{\a+\b+3}(-1)^{\ell}
\frac{(-\a-\b-3)_{\ell}}{\ell!}\l(\frac{x-1}{2}\r)^{2\a+2\b+6-\ell}\nn
&=&\frac{\a+\b+2}{(\a+1)(\b+1)}
\frac{2^{2\a+2\b+6}}{(\a+\b+1)!\,(\a+\b+3)!}(1-x^2)\nn
& &{}\hspace{5cm}\times
\l(\frac{x-1}{2}\r)^{\a+\b+3}\l(\frac{x+1}{2}\r)^{\a+\b+3}\nn
&=&-\frac{\a+\b+2}{(\a+1)(\b+1)}
\frac{(x^2-1)^{\a+\b+4}}{(\a+\b+1)!\,(\a+\b+3)!}.\n\eea

Now we use the fact that
$$C_i(x)=(x^2-1)c_i(x)+2(i+1)xc_{i+1}(x)+(i+1)(i+2)c_{i+2}(x),
\;i=2,3,4,\ldots$$
to conclude that
$$C_{2\a+2\b+6}(x)=(x^2-1)c_{2\a+2\b+6}(x),$$
which leads to (\ref{ckop}).

\section{Some remarks}

Let $\a>-1$ and $\b>-1$.

The coefficients $\l\{a_i(x)\r\}_{i=1}^{\infty}$ and
$\l\{b_i(x)\r\}_{i=1}^{\infty}$ can also be computed in the same way as we
computed the coefficients $\l\{c_i(x)\r\}_{i=1}^{\infty}$. Consider the system
of equations $S_4=0$. First we use (\ref{diff}) to find from (\ref{defR1})
that
$$\R=\frac{(\a+2)_{n-1}(\a+\b+2)_n}{2(\b+1)_{n-1}\,n!}(x-1)
P_{n-1}^{(\a+2,\b)}(x),\;n=1,2,3,\ldots.$$
Now we use the fact that
$$\frac{(\a+2)_{n-1}(\a+\b+2)_n}{2(\b+1)_{n-1}\,n!}\ne 0,
\;n=1,2,3,\ldots$$
to conclude that
$$\sum_{i=0}^{\infty}b_i(x)D^i\l[(x-1)P_{n-1}^{(\a+2,\b)}(x)\r]=0,
\;n=1,2,3,\ldots.$$
Hence, by shifting $n$ we obtain
$$\sum_{i=1}^{\infty}b_i(x)D^i\l[(x-1)P_n^{(\a+2,\b)}(x)\r]=
-b_0(n+1,\a,\b)(x-1)P_n^{(\a+2,\b)}(x),\;\ndots.$$
Note that for $\ndots$ we have
$$D^i\l[(x-1)P_n^{(\a+2,\b)}(x)\r]=
(x-1)D^iP_n^{(\a+2,\b)}(x)+iD^{i-1}P_n^{(\a+2,\b)}(x),\;i=1,2,3,\ldots.$$
Hence we obtain
\be\la{sys2}\sum_{i=0}^{\infty}B_i(x)D^iP_n^{(\a+2,\b)}(x)=
-b_0(n+1,\a,\b)(x-1)P_n^{(\a+2,\b)}(x),\;\ndots,\ee
where
$$B_i(x)=\l\{\ba{ll}b_1(x), &i=0\\[5mm]
(x-1)b_i(x)+(i+1)b_{i+1}(x), &i=1,2,3,\ldots.\ea\r.$$
Note that the system of equations (\ref{sys2}) has the form (\ref{inv2}). So
we may apply (\ref{opl2}) and use (\ref{bnul}) to conclude that for
$i=0,1,2,\ldots$ we have
\bea B_i(x)&=&-(x-1)2^i\sum_{j=0}^i\frac{\a+\b+2j+3}{(\a+\b+j+3)_{i+1}}
b_0(j+1,\a,\b)P_{i-j}^{(-\a-i-3,-\b-i-1)}(x)P_j^{(\a+2,\b)}(x)\nn
&=&\la{B}-(\a+\b+2)(x-1)2^i\nn
& &{}\hspace{1cm}\times\sum_{j=0}^i\frac{\a+\b+2j+3}{(\a+\b+j+3)_{i+1}}
\frac{(\a+3)_j(\a+\b+3)_j}{(\b+1)_j\,j!}\nn
& &{}\hspace{5cm}\times P_{i-j}^{(-\a-i-3,-\b-i-1)}(x)P_j^{(\a+2,\b)}(x).
\eea

As before we can deduce that for $i=0,1,2,\ldots$
\bea B_i(x)&=&-(\a+\b+2)(x-1)(-2)^i\nn
& &{}\hspace{1cm}\times\sum_{\ell=0}^i\frac{(\a+3)_{i-\ell}(-\a-2)_{i-\ell}}
{(\b+1)_{i-\ell}(i-\ell)!\,\ell!\,(i-\ell)!}\nn
& &{}\hspace{3cm}\times
\hyp{3}{2}{-\ell,\a+i-\ell+3,\a+\b+3}{\b+i-\ell+1,i-\ell+1}{1}
\l(\frac{x-1}{2}\r)^{\ell}.\n\eea

Now we use (\ref{A1}) and (\ref{B1}) with $z=x-1$ to find that
\bea b_i(x)&=&\frac{1}{i!}\sum_{j=0}^{i-1}
(-1)^{i-j-1}j!\,(x-1)^{i-j-1}B_j(x)\nn
&=&\frac{1}{i!}\sum_{j=0}^{i-1}(-1)^j(i-j-1)!\,(x-1)^jB_{i-j-1}(x),
\;i=1,2,3,\ldots\n\eea
which leads to (\ref{b}) after changing the order of summations and using
the Vandermonde summation formula (\ref{Van}) as before.

In a similar way the coefficients $\l\{a_i(x)\r\}_{i=1}^{\infty}$ can be
computed from the system of equations $S_2=0$. In that case we would need
(\ref{A1}) and (\ref{B1}) with $z=x+1$, but it is easier to use the symmetry
relation (\ref{symab}) of course.

In \cite{Bav2} H.~Bavinck found the following interesting formula involving
Laguerre polynomials~:
$$\sum_{k=j}^ik^sL_{i-k}^{(-\a-i-1)}(-x)L_{k-j}^{(\a+j)}(x)
=(-x)^s\delta_{i,j+2s},\;i\ge j+2s,\;i,j,s\in\{0,1,2,\ldots\},$$
which holds for all $\a$. In \cite{Bav4} he found an analogue of this
formula involving Jacobi polynomials~:
\bea & &2^{i-j}\sum_{k=j}^i
\frac{\a+\b+2k+1}{(\a+\b+j+k+1)_{i-j+1}}\l[k(k+\a+\b+1)\r]^s
P_{i-k}^{(-\a-i-1,-\b-i-1)}(x)P_{k-j}^{(\a+j,\b+j)}(x)\nn
& &{}\hspace{3cm}
=(x^2-1)^s\delta_{i,j+2s},\;i\ge j+2s,\;i,j,s\in\{0,1,2,\ldots\},\n\eea
which holds for $-(\a+\b+2)\notin\{0,1,2,\ldots\}$. The case $\a+\b+1=0$
must be understood by continuity. This formula can be applied to (\ref{C1})
and (\ref{B}). Since we have
$$\frac{(\a+3)_j(\a+\b+3)_j}{(\b+1)_j\,j!}=
\frac{(j+1)_{\a+2}(\b+j+1)_{\a+2}}{(\b+1)_{\a+2}(\a+2)!},
\;\a\in\{0,1,2,\ldots\}$$
and
$$(j+1)_{\a+2}(\b+j+1)_{\a+2}=\prod_{k=1}^{\a+2}
\l[j(j+\a+\b+3)+k(\a+\b+3-k)\r],\;\a\in\{0,1,2,\ldots\}$$
this implies that for $\a\in\{0,1,2,\ldots\}$
$$B_i(x)=-(\a+\b+2)(x-1)\frac{(x^2-1)^{\a+2}}{(\b+1)_{\a+2}(\a+2)!}
\delta_{i,2\a+4},\;i\ge 2\a+4.$$
Hence, for $\a\in\{0,1,2,\ldots\}$ we have
$$B_i(x)=0,\;i>2\a+4\;\textrm{ and }\;B_{2\a+4}(x)=-(x-1)
\frac{(x^2-1)^{\a+2}}{(\b+1)_{\a+1}(\a+2)!},$$
which leads to (\ref{bkop}) eventually. In a similar way we find for
$\a,\b\in\{0,1,2,\ldots\}$
$$\frac{(\a+\b+4)_{j+1}(\a+\b+4)_j}{(j+1)!\,j!}=
\frac{(j+2)_{\a+\b+3}(j+1)_{\a+\b+3}}{(\a+\b+3)!\,(\a+\b+3)!}$$
and
$$(j+2)_{\a+\b+3}(j+1)_{\a+\b+3}=\prod_{k=1}^{\a+\b+3}
\l[j(j+\a+\b+5)+k(\a+\b+5-k)\r],$$
which implies that in view of (\ref{C1}) we have
$$C_i(x)=-\frac{(\a+\b+2)^2(\a+\b+3)}{(\a+1)(\b+1)}
\frac{(x^2-1)^{\a+\b+4}}{(\a+\b+3)!\,(\a+\b+3)!}\delta_{i,2\a+2\b+6},
\;i\ge 2\a+2\b+6.$$
Hence, for $\a,\b\in\{0,1,2,\ldots\}$ we have
$$C_i(x)=0,\;i>2\a+2\b+6\;\textrm{ and }\;C_{2\a+2\b+6}(x)=
-\frac{\a+\b+2}{(\a+1)(\b+1)}
\frac{(x^2-1)^{\a+\b+4}}{(\a+\b+1)!\,(\a+\b+3)!}$$
as before.

If we set $M=0$ into (\ref{DV}) we get the differential equation
\bea\la{Dvlim}& &N\sum_{i=0}^{\infty}b_i(x)y^{(i)}(x)+(1-x^2)y''(x)\nn
& &{}\hspace{1cm}{}+\l[\b-\a-(\a+\b+2)x\r]y'(x)
+n(n+\a+\b+1)y(x)=0,\eea
satisfied by the polynomials $\set{P_n^{\a,\b,0,N}(x)}$.
From the limit relation (\ref{lim}) it follows that
$$\lim_{\b\rightarrow\infty}\frac{(-2)^i}{\b^i}\l(D^i
P_n^{\a,\b,0,N}\r)\l(1-\frac{2x}{\b}\r)=D^iL_n^{\a,N}(x),
\;\ndots,\;i=0,1,2,\ldots,$$
where $L_n^{\a,N}(x)$ denotes the generalized Laguerre polynomial
considered in \cite{Dvlag}. Note that from (\ref{bnul}) we easily
find that
$$\lim_{\b\rightarrow\infty}\frac{b_0(n,\a,\b)}{\b}=
\l({n+\a+1 \atop n-1}\r),\;\ndots.$$
Now we use the Vandermonde summation formula (\ref{Van}) to find
from (\ref{b})
\bea & &\lim_{\b\rightarrow\infty}\frac{\b^{i-1}}{(-2)^i}
b_i(\a,\b,1-2x/\b)\nn
&=&\sum_{\ell=0}^{i-1}\frac{(\a+3)_{i-\ell-1}(-\a-2)_{i-\ell-1}}
{(i-\ell)!\,(i-\ell-1)!\,\ell!}\,
\hyp{2}{1}{-\ell,\a+i-\ell+2}{i-\ell+1}{1}(-x)^{\ell+1}\nn
&=&\sum_{\ell=0}^{i-1}
\frac{(\a+3)_{i-\ell-1}(-\a-2)_{i-\ell-1}(-\a-1)_{\ell}}
{i!\,(i-\ell-1)!\,\ell!}(-x)^{\ell+1}\nn
&=&\frac{1}{i!}\sum_{j=1}^i(-1)^{i+j+1}\l({\a+1 \atop j-1}\r)
\l({\a+2 \atop i-j}\r)(\a+3)_{i-j}x^j,\;i=1,2,3,\ldots.\n\eea
Hence, if we set $y(x)=P_n^{\a,\b,0,N}(x)$ into the differential
equation (\ref{Dvlim}), change $x$ by $1-2x/\b$, divide by $\b$
and take the limit $\b\rightarrow\infty$ we obtain the differential
equation for the polynomials $\set{L_n^{\a,N}(x)}$ which was found
in \cite{Dvlag}.

In \cite{Madrid} and \cite{Symjac} we found all differential
equations of the form (\ref{DVSym}) satisfied by the polynomials
$\set{\SGP}$, where $\a>-1$ and $M\ge 0$. We emphasize that these
differential equations are not of the form (\ref{DV}).
The differential equation (\ref{DV}) leads to another one after
setting $\b=\a$ and $N=M$.

In \cite{Symjac} we also found differential equations for the polynomials
$\set{P_n^{\a,\pm\frac{1}{2},0,N}}$, where $\a>-1$ and $N\ge 0$. It can be
shown that these do coincide with (\ref{DV}) after setting $M=0$ and
$\b=\pm\frac{1}{2}$. For $\a\in\{0,1,2,\ldots\}$ and $N>0$ these
differential equations have finite order $2\a+4$.

Finally we can correct a table conjectured in \cite{Everitt} listing the
cases for which the polynomials $\set{\GP}$ satisfy a finite order
differential equation of the form (\ref{DV}) with minimal order~:

\vspace{5mm}

$$\begin{array}{|c|c|c|}
\hline
M,N & \a,\b & \textrm{order} \\
\hline
\hline
M=0,\;N=0 & \a>-1,\;\b>-1 & 2\\
\hline
N=M>0 & \b=\a\in\{0,1,2,\ldots\} & 2\a+4\\
\hline
M=0,\;N>0 & \a\in\{0,1,2,\ldots\},\;\b>-1 & 2\a+4\\
\hline
M>0,\;N=0 & \a>-1,\;\b\in\{0,1,2,\ldots\} & 2\b+4\\
\hline
M>0,\;N>0 & \a,\b\in\{0,1,2,\ldots\} & 2\a+2\b+6\\
\hline
\end{array}$$

\vspace{1cm}

Menelaoslaan 4, 5631 LN Eindhoven, The Netherlands

\vspace{5mm}

Delft University of Technology, Faculty of Information Technology and
Sytems,

P.O. Box 5031, 2600 GA Delft, The Netherlands, e-mail :
koekoek@twi.tudelft.nl


\begin{thebibliography}{99}
\bibitem{Bailey} {\sc W.N. Bailey :} {\em Generalized hypergeometric
series.} Cambridge University Press, Cambridge, 1935.
\bibitem{Bav1} {\sc H. Bavinck :} {\em A direct approach to Koekoek's
differential equation for generalized Laguerre polynomials.} Acta
Mathematica Hungarica {\bf 66}, 1995, 247-253.
\bibitem{Bav2} {\sc H. Bavinck :} {\em A new result for Laguerre
polynomials.} Journal of Physics A {\bf 29}, 1996, L277-L279.
\bibitem{Bav3} {\sc H. Bavinck :} {\em Differential and difference
operators having orthogonal polynomials with two linear perturbations as
eigenfunctions.} Journal of Computational and Applied Mathematics {\bf 92},
1998, 85-95.
\bibitem{Bav4} {\sc H. Bavinck :} {\em Differential operators
having Sobolev-type Gegenbauer polynomials as eigenfunctions.}
To appear.
\bibitem{Char} {\sc H. Bavinck \& R. Koekoek :} {\em On a difference
equation for generalizations of Charlier polynomials.} Journal of
Approximation Theory {\bf 81}, 1995, 195-206.
\bibitem{Chihara} {\sc T.S. Chihara :} {\em An introduction to orthogonal
polynomials.} Mathematics and Its Applications {\bf 13}, Gordon and Breach,
New York, 1978.
\bibitem{Everitt} {\sc W.N. Everitt \& L.L. Littlejohn :} {\em Orthogonal
polynomials and spectral theory : a survey.} Orthogonal Polynomials and
their Applications, Volume~{\bf 9} of IMACS Annals on Computing and Applied
Mathematics (Editors~: C. Brezinski, L. Gori \& A. Ronveaux).
J.C.~Baltzer~AG, Basel, 1991, 21-55.
\bibitem{Dvlag} {\sc J. Koekoek \& R. Koekoek :} {\em On a differential
equation for Koornwinder's generalized Laguerre polynomials.} Proceedings of
the American Mathematical Society {\bf 112}, 1991, 1045-1054.
\bibitem{Madrid} {\sc J. Koekoek \& R. Koekoek :} {\em Finding differential
equations for symmetric generalized ultraspherical polynomials by using
inversion methods.} Proceedings of the International Workshop on Orthogonal
Polynomials in Mathematical Physics (Legan\'es, 1996), Universidad Carlos
III de Madrid, Legan\'es, 1997, 103-111.
\bibitem{Invjac} {\sc J. Koekoek \& R. Koekoek :} {\em The Jacobi inversion
formula.} Complex Variables, to appear.
\bibitem{Soblag} {\sc J. Koekoek, R. Koekoek \& H. Bavinck :} {\em On
differential equations for Sobolev-type Laguerre polynomials.} Transactions
of the American Mathematical Society {\bf 350}, 1998, 347-393.
\bibitem{Search} {\sc R. Koekoek :} {\em The search for differential
equations for certain sets of orthogonal polynomials.} Journal of
Computational and Applied Mathematics {\bf 49}, 1993, 111-119.
\bibitem{Symjac} {\sc R. Koekoek :} {\em Differential equations for
symmetric generalized ultraspherical polynomials.} Transactions of the
American Mathematical Society {\bf 345}, 1994, 47-72.
\bibitem{Askey} {\sc R. Koekoek \& R.F. Swarttouw :} {\em The Askey-scheme
of hypergeometric orthogonal polynomials and its $q$-analogue.} Delft
University of Technology, Faculty of Technical Mathematics and Informatics,
report no. {\bf 98-17}, 1998.
\bibitem{Koorn} {\sc T.H. Koornwinder :} {\em Orthogonal polynomials with
weight function $(1-x)^{\a}(1+x)^{\b}+M\delta(x+1)+N\delta(x-1)$.}
Canadian Mathematical Bulletin {\bf 27}(2), 1984, 205-214.
\bibitem{A.M.Krall} {\sc A.M. Krall :} {\em Orthogonal polynomials
satisfying fourth order differential equations.} Proceedings of the Royal
Society of Edinburgh {\bf 87A}, 1981, 271-288.
\bibitem{H.L.Krall1} {\sc H.L. Krall :} {\em Certain differential equations
for Tchebycheff polynomials.} Duke Mathematical Journal {\bf 4}, 1938,
705-718.
\bibitem{H.L.Krall2} {\sc H.L. Krall :} {\em On orthogonal polynomials
satisfying a certain fourth order differential equation.} The Pennsylvania
State College Studies, No. 6, 1940.
\bibitem{Littlejohn} {\sc L.L. Littlejohn :} {\em The Krall polynomials~: A
new class of orthogonal polynomials.} Quaestiones Mathematicae {\bf 5},
1982, 255-265.
\bibitem{Luke} {\sc Y.L. Luke :} {\em The special functions and their
approximations II.} Academic Press, New York, 1969.
\bibitem{Slater} {\sc L.J. Slater :} {\em Generalized hypergeometric
functions.} Cambridge University Press, Cambridge, 1966.
\bibitem{Szego} {\sc G. Szeg{\"o} :} {\em Orthogonal polynomials.}
American Mathematical Society Colloquium Publications {\bf 23} (1939),
Fourth edition, Providence, Rhode Island, 1975.
\end{thebibliography}
\end{document}